\documentclass[a4paper, english, 12pt, twoside]{article}
\usepackage[a4paper, total={7in, 10in}, twoside, bindingoffset=6mm]{geometry}

\usepackage[utf8]{inputenc}

\usepackage{babel}
\usepackage{hyphsubst}
\HyphSubstIfExists{english-x-latest}{
  \HyphSubstLet{english}{english-x-latest}}{}
\usepackage{lmodern}

\usepackage{textcomp}
\usepackage{microtype}

\usepackage{setspace}

\usepackage{blindtext}
\usepackage{todonotes}

\usepackage{authblk}

\usepackage{tikz}
\usepackage{neuralnetwork}
\usepackage{pgfplots}
\pgfplotsset{compat=1.17}
\usetikzlibrary{positioning, arrows.meta}

\usepackage{overpic}
\usepackage[usestackEOL]{stackengine}

\usepackage{tabularx}

\usepackage{verbatim}

\usepackage[shortlabels]{enumitem}

\usepackage{bm}

\usepackage{graphicx}
\graphicspath{{./}}
\usepackage{subcaption}
\usepackage{svg}
\svgsetup{inkscape=overwrite,inkscapeexe="C:/Program Files/Inkscape/bin/inkscape.exe",inkscapepath=svgsubdir}
\svgpath{{./}}
\usepackage[font=scriptsize,labelfont=bf]{caption}
\usepackage{float}
\usepackage{lipsum}

\usepackage{amsmath}
\usepackage{amsfonts}
\usepackage{amsthm}
\usepackage{amssymb}

\allowdisplaybreaks

\usepackage[per-mode = symbol]{siunitx}
\DeclareSIUnit[quantity-product = \,]
    {\darcy}{\text{D}}
\DeclareSIUnit[quantity-product = \,]
    {\bar}{\text{bar}}

\usepackage{mathtools}
\mathtoolsset{centercolon}

\usepackage{rotating}

\usepackage[nottoc,numbib]{tocbibind}

\newcommand{\R}{\mathbb{R}}

\DeclareMathOperator*{\argmin}{arg\,min}

\theoremstyle{definition}

\theoremstyle{plain}
\newtheorem*{theorem*}{Theorem} 
\newtheorem*{corollary*}{Corollary}
\newtheorem*{lemma*}{Lemma}
\newtheorem*{proposition*}{Proposition}

\theoremstyle{definition} 
\newtheorem*{definition*}{Definition}

\usepackage{xcolor, soul}
\sethlcolor{red}
\definecolor{blue}{rgb}{240,80,30}

\usepackage{color}
\usepackage{listings}

\usepackage{csquotes}

\usepackage[round,authoryear]{natbib}
\usepackage{etoolbox}

\makeatletter
\patchcmd{\blx@imc@online}{\blx@sf@dot}{\blx@sf@none}{}{}
\makeatother

\newcounter{nalg}[section]
\renewcommand{\thenalg}{\thesection .\arabic{nalg}}
\DeclareCaptionLabelFormat{algocaption}{Algorithm \thenalg}

\lstnewenvironment{algorithm}[1][]
{   
    \refstepcounter{nalg}
    \captionsetup{labelformat=algocaption,labelsep=colon,justification=raggedright,singlelinecheck=off}
rr  \lstset{
        mathescape=true,
        frame=tB,
        numbers=left, 
        numberstyle=\tiny,
        basicstyle=\scriptsize, 
        keywordstyle=\color{black}\bfseries,
        keywords={,input, output, return, datatype, function, in, if, else, for, each, while, begin, end, sum, abs,initialize,generate, elif,}
        numbers=left,
        xleftmargin=.04\textwidth,
        #1
	}
}
{}

\usepackage[
	bookmarksopen,
	bookmarksnumbered,
	pdftitle={A machine-learned near-well model in OPM Flow},
	pdfauthor={Peter von Schultzendorff},
	pdfkeywords={CO2 storage, reservoir simulator, neural networks, OPM Flow, OPM, hybrid models, automatic differentiation, near-well model},
	pdfstartview={FitH},
]
{hyperref}

\usepackage{doi}

\usepackage[noabbrev,capitalize]{cleveref}

\usepackage{subfiles}

\title{A machine-learned near-well model in OPM Flow}

\author[1]{Peter von Schultzendorff}
\author[2]{Tor Harald Sandve}
\author[2]{Birane Kane}
\author[2]{David Landa-Marbán}
\author[1]{Jakub Wiktor Both}
\author[1]{Jan Martin Nordbotten}
\affil[1]{Universitetet i Bergen, Norway}
\affil[2]{NORCE, Norway}

\begin{document}

	\maketitle

	\begin{abstract}
		Recently, a focal point of research in reservoir simulation and other fields have been hybrid approaches that couple classical physics-based techniques with data-driven models. The latter offer high fidelity to real-world phenomena and fast inference, which enables both faster simulations and the possibility to model not fully understood phenomena to good precision. To facilitate hybrid models to the maximum extent, reservoir simulators require full integration of machine-learning components into the numerical machinery.

		Furthermore, comprehensive automatic differentiation of simulators has received much interest in the past years. This enables, e.g., efficient and flexible strategies for nonlinear solvers, inverse problems, and optimization problems. To retain these workflows with hybrid models, it is of particular importance to combine the automatic differentiation functionalities of employed machine-learning frameworks with those of the simulators. 

		Here, we present the first integration of neural networks into the reservoir simulator OPM Flow. OPM Flow as a competitive, open-source, high performance, and production-ready simulator sets high standards for such a task. We devise a solution that integrates seamlessly into existing workflows, is simple from the user perspective, yet adaptable to different use cases. Models are trained in TensorFlow, stored, and then loaded into OPM, where they can be simply accessed as functions, fully integrated into the automatic differentiation framework of OPM. The result is an efficient and flexible framework for hybrid modeling, which ensures full automatic differentiation capabilities for the models.

		To showcase possible applications of this new framework, we introduce a novel, data-driven near-well model. In reservoir simulations, it is common practice to employ near-well models to accurately capture the singular character of the pressure gradients in proximity to wells. Commonly used are the Peaceman well model, its various extensions, or local grid refinement around the wells. These techniques have different limitations, as the former approaches are only accurate in simplified flow regimes, while the latter can be computationally expensive. Our approach addresses these challenges by integrating a machine-learned near-well model into OPM Flow via the aforementioned framework. The key idea is to obtain a data-driven inflow-performance relation by conducting fine-scale ensemble simulations of the near-well region. A neural network is then trained on the data and infers the inflow-performance relation at simulation time. Tested on relevant examples for CO\textsubscript{2} storage, the method offers high fidelity to fine-scale results at low computational overload.

		This application demonstrates the potential that the OPM Flow - Neural Network framework offers for hybrid modeling.
	\end{abstract}

	\section{Introduction}\label{section:introduction}
		
	In recent years, a focal point of research in reservoir simulation has been the facilitation of machine-learning (ML) models. By learning from field data, experiments, and/or previous simulations, they enable both faster computation and higher accuracy. ML models can be facilitated at various parts of the reservoir simulation workflow. Perhaps most importantly, they can improve individual simulation runs, either by reducing computational complexity or by providing high fidelity to real-world data. Furthermore, ML is used in history matching and to improve production forecast/optimization models \citep{Samnioti:ApplicationsMachineLearning2023,Samnioti:ApplicationsMachineLearning2023a}. In all applications, ML is a valuable asset because it allows for accurate and computationally efficient representation of physical processes, under the assumption that these processes are expressed in the underlying data. In particular when modeling unknown and/or not-fully-understood physics that are absent in current models, this proves useful. The data-driven approach removes the necessity to describe such physical processes as first principle laws, while retaining high accuracy.

	When enhancing single simulation runs with ML, the main goals are to reduce computational complexity and to improve accuracy. There exist a variety of approaches to achieve these goals. Surrogate reservoir models mimic the behavior of full-field models simulation with some ML approach. These models infer either pressure and saturation fields over the full domain \citep{Shahkarami:ModelingPressureSaturation2014,Bao:DevelopmentProxyModels2018,Kiaerr:EvaluationDataDrivenFlow2020}, or singular simulation values such as well pressure \citep{Memon:SurrogateReservoirModelingprediction2014}. ML models can also be combined with classical physics-based models in hybrid models. For example, they can be integrated into existing simulators, where they replace only certain parts of the simulation instead of providing the full simulations results themselves. Notable applications of this are the solution of stability and phase split problems in compositional models with ML \citep{Li:AccelerationNVTFlash2019}, ML models that provide PVT properties in black oil models \citep{Varotsis:NovelNonIterativeMethod1999}, and data-driven relative permeability models \citep{Arigbe:RealtimeRelativePermeability2019}. Another utilization of ML models is to improve numerical solvers. One exemplary approach is to provide better initial guesses for solvers to reduce the number of iterations during solving \citep{Lechevallier:HybridNewtonMethod2023}. Notably, the final simulation result is still obtained from the classic physics-based simulator.

	Another class of hybrid approaches that facilitate ML are digital twins, i.e., full digital models of physical assets that continuously integrate real-time data to update the model. For large systems, the real-time requirement necessitates the combination of ML and physics-based modeling \citep{Rasheed:DigitalTwinValues2020}. An example is the \emph{PoroTwin} project \citep{Keilegavlen:PoroTwinDigitalTwin2023}, which corrects simulation results with ML to be closer to the real-world model. Alternatively, ML can modify a model's solution to integrate complimentary physics not included in said model. An example is the integration of reservoir souring via ML into the surrogate \emph{FlowNet} model \citep{Leeuwenburgh:HybridDataPhysicsFramework2024}.

\begin{figure}[h!]
    	\color{gray}
    	\centering
    	\stackinset{l}{121pt}{b}{287pt}{(\cite{preCICEv2})}{
    	\stackinset{l}{133pt}{b}{104pt}{(\cite{Leeuwenburgh:HybridDataPhysicsFramework2024})}{
    	\stackinset{l}{164pt}{b}{381pt}{(\cite{Shahkarami:ModelingPressureSaturation2014})}{
    	\stackinset{l}{176pt}{b}{149pt}{(\cite{Bao:DevelopmentProxyModels2018})}{
    	\stackinset{l}{105pt}{b}{225pt}{(\cite{Keilegavlen:PoroTwinDigitalTwin2023})}{
    	\stackinset{l}{198pt}{b}{57pt}{(\cite{Lechevallier:HybridNewtonMethod2023})}{
		\stackinset{l}{323pt}{b}{100pt}{(\cite{Despres:MachineLearningDesign2020})}{\includeinkscape[inkscapelatex=false, width=.75\linewidth]{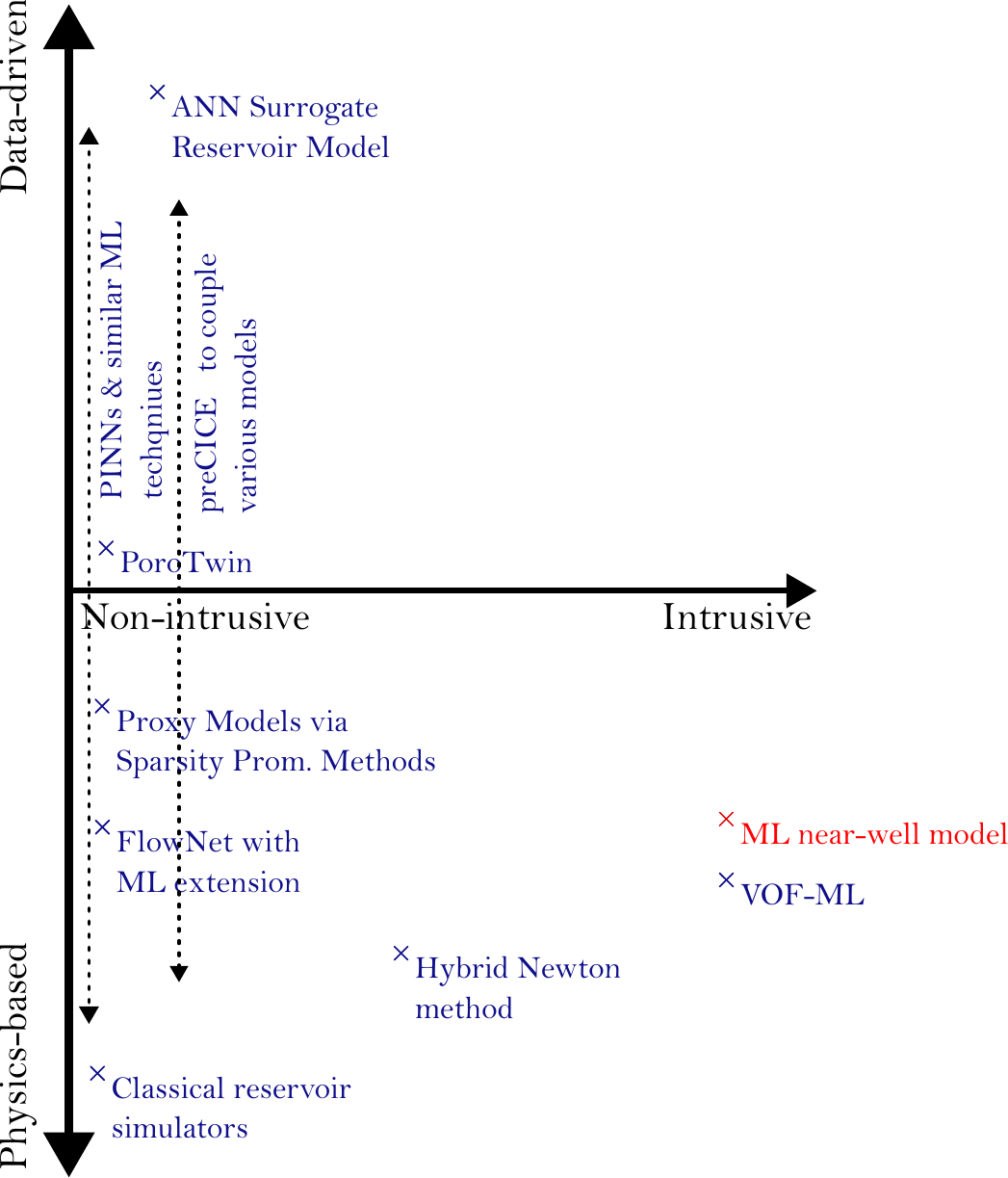}}}}}}}}
        \caption[Methods]{A brief overview of existing physics-based, data-driven, and hybrid models and simulation frameworks. On the \(y\)-axis models are arranged by whether they are predominantly physics- or data-driven, while on the \(x\)-axis they are arranged by how intrusive their implementation is.}
        \label{figure:methods}
	\end{figure}

	Hybrid approaches that make use of both physics-based simulators and ML models can be categorized by the intrusiveness of their implementation. In an intrusive approach, the ML components and the classical model are tightly integrated, i.e., either the ML components are (partly) integrated into the physics-based simulator or vice versa. In a non-intrusive approach, both models operate on their own and communicate only via an interface. \cref{figure:methods} displays physics-based, data-driven, and hybrid models in reservoir simulations arranged by the intrusiveness of their implementation. An interface approach retains a high flexibility, as models on either side of the interface can be exchanged easily. An example for such a framework is the \emph{preCICE} library that can couple together various existing simulators via adapters to enable complex multi-physics coupling \citep{preCICEv2}. However, in particular for hybrid models that couple physics-based and ML models, e.g, \citep{Li:AccelerationNVTFlash2019,Varotsis:NovelNonIterativeMethod1999,Arigbe:RealtimeRelativePermeability2019}, an intrusive approach can be of advantage. This is, because overhead for communication between different model components is reduced. Additionally, in recent years, a focus has been laid on comprehensive automatic differentiation capabilities of simulators that improve flexibility and extensibility to new models and numerical methods \citep{Rasmussen:OpenPorousMedia2021}. Intrusive implementations can facilitate full automatic differentiation capacities for hybrid models. However, much of the work on ML components just showcases a proof-of-concept ML model itself without integration into a simulator \citep{Li:AccelerationNVTFlash2019,Varotsis:NovelNonIterativeMethod1999,Arigbe:RealtimeRelativePermeability2019}. This highlights the challenges of integrating ML models into reservoir simulators and the need for a simple and user-friendly framework for hybrid models.

	Motivated by these demands, we present the first integration of artificial neural networks (ANNs; short NNs) into the reservoir simulator \emph{Open Porous Media Flow} (OPM Flow) \citep{Rasmussen:OpenPorousMedia2021}. OPM Flow as a competitive, open-source, high performance, and production-ready simulator sets high standards for such a task. We devise the \emph{OPM Flow-NN framework}, which integrates seamlessly into existing workflows, is simple from the user perspective, yet adaptable to different use cases. Models are trained offline in \emph{Keras} \citep{Chollet:Keras2015}, stored, and then loaded into OPM Flow, where they can be simply accessed as native functions.

	Embedding NNs trained in Keras/Python into a C++ application works by saving the model in an appropriate format and using an automated tool to generate corresponding C++ code. The generated C++ code is then further customized and incorporated into OPM Flow. We build our work on the tools \emph{Kerasify} \citep{RobertW.RosePaulMaevskikh:Kerasify2017} and \emph{Keras2cpp} \citep{RobertW.RoseGeorgyPerevozchikov:Keras2cpp2019}, but make customizations that enable integration into OPM Flow. We point out that this procedure enables the NNs to be fully integrated into the automatic differentiation framework of OPM Flow. The result is an efficient and flexible framework for hybrid modeling, which ensures comprehensive automatic differentiation capabilities.

	As an application of the OPM Flow-NN framework, we showcase a novel, machine-learned near-well model. In reservoir simulations, pressure gradients attain a singular character in proximity of production and injection wells. These gradients are not adequately captured by typical grid sizes for simulations, which are in the range of tenths to hundreds of meters. Accurate representation of near-well phenomena has therefore been a focal point of research over the past decades. Typically, near-well models are employed that allow for accurate modeling of the well flow rate \(q\) and the well pressure \(p_{well}.\) A common idea in many near-well models is the introduction of a well index \(WI,\) which defines the transmissibility between a well connection \(i\) and the grid block it is encompassed by. The well index linearly relates \(p_{well},\) \(q,\) and the well block pressure \(p_{i},\) here exemplary displayed for an injection well
	\begin{subequations}\label{equation:peaceman-original-model}
	\begin{equation}
    	q = WI (p_{well} - p_i). \label{equation:peaceman-original-model:linear-relation}
	\end{equation}
    One of the most known and classical derivations of the well index is the Peaceman well model \citep{Peaceman:InterpretationWellBlockPressures1978}. It provides, under the assumption of single-phase, incompressible, radial flow in a homogeneous and isotropic reservoir, an analytic and thus exact expression for the well index
	\begin{equation}\label{equation:peaceman-original-model:WI}
    	WI = \frac{2 \pi k h \frac{\rho}{\mu}}{\ln\frac{r_e}{r_w}}.
	\end{equation}
	\end{subequations}
	Here, \(k\) is (isotropic) permeability, \(h\) is the height of the grid block containing the well section of interest, \(\rho\) and \(\mu\) are phase density and viscosity, respectively, and \(r_w\) is the well radius. The term \(r_e\) denotes the \emph{equivalent well radius,} i.e., the distance from the well where the exact pressure value equals the grid block pressure. For flow on a repeated five-spot pattern on a quadratic grid, \(r_e\) can be computed exactly in terms of the cell size \(\Delta x,\) hence \(WI\) \eqref{equation:peaceman-original-model:WI} gives an exact relation between well and block pressure \citep{Peaceman:InterpretationWellBlockPressures1978}.

	Peaceman's original model has been extended to a variety of different cases. These include specific models for non-square blocks and anisotropic reservoirs \citep{Peaceman:InterpretationWellBlockPressures1983}, off-center wells and multiple wells per block \citep{Peaceman:InterpretationWellblockPressures1990}, slanted wells \citep{Aavatsmark:WellIndexReservoir2003}, and horizontal wells \citep{Peaceman:InterpretationWellBlockPressures1983, Cao:ModelingSemisteadyState2018}. While these (semi-)analytic approaches provide exact results under simplified flow regimes, they also entail several drawbacks. We highlight in particular the loss of accuracy in the near-well region for an increasing degree of heterogeneity, as well as the need for case-specific analysis for different settings.

	Flow-based numerical upscaling methods retain the well index approach \eqref{equation:peaceman-original-model:linear-relation}, but address reservoirs with significant heterogeneities. The flow and pressure values from numerical fine-scale near-well solutions to the pressure equation are averaged to obtain a representation of the heterogeneity (in terms of upscaled values for transmissibility and well index) at coarse-scale \citep{Ding:ScalingupVicinityWells1995, Wolfsteiner:NearWellRadialUpscaling2002}. While upscaling methods resolve heterogeneities, they are inaccurate in cases such as transient flow or in the presence of multiple phases near the wellbore. In the latter, the flow-performance relation is largely dependent on the phase mobilities, which makes fine resolution of the phase saturation fields essential for accurate modeling. Multiphase upscaling mitigates this issue \citep{Nakashima:AccurateRepresentationNearwell2010}, but is limited in expressing complex relations due to its conception as a lookup-table approach. Alternatively, transient effects can be included via a time-dependent well index. However, the right choice of the equivalent well radius in such cases is not straightforward \citep{Al-Mohannadi:GridSystemRequirementsNumerical2007}.

	Local grid refinement (LGR) replaces the well index approach \eqref{equation:peaceman-original-model:linear-relation} by resolving near-well regions with sufficiently fine grid cells during full simulations \citep{Karimi-Fard:AccurateResolutionNearWell2012}. The method is accurate in general cases, but may require a large number of additional grid cells and thus become computationally expensive.

	The presented machine-learned near-well model retains the accuracy and flexibility of LGR, while shifting the computational load from \emph{online} (during simulation) to \emph{offline} (before the simulation). The key idea is to replace the well index \(WI\) \eqref{equation:peaceman-original-model:linear-relation} from the (semi-)analytic and upscaling methods by a data-driven expression. With this, a machine-learning routine takes the place of case-specific analytic derivation, moving the difficulty of obtaining an accurate and robust model to the availability of data and the training. Ensemble simulations of the near-well region in different flow regimes are conducted and upscaled to form a training dataset. Afterward, a NN is trained on the data to learn an improved \(WI.\) It is integrated into OPM Flow with the OPM Flow-NN framework. The resulting workflow for near-well modeling is accurate across different flow regimes, robust to variations in model parameters, flexible in terms of potential scenarios, and fast at inference time.

	This work proceeds as follows. First, we describe the integration of NNs into OPM Flow (section "\nameref{section:OPM Flow-NN}"). Second, we introduce the data-driven near-well model (section "\nameref{section:near well model}"). To test the model, we consider one single-phase and two two-phase injection scenarios, the first of the latter with strong transient effects affecting the multiphase flow regime and the second with vertical flow between heterogeneous layers. These results are discussed (section "\nameref{section:simulations}"), before concluding with remarks and outlook (section "\nameref{section:conclusion}").

	\section{OPM Flow-NN framework}\label{section:OPM Flow-NN}
		
	The OPM Flow-NN framework facilitates simple and flexible use of NNs in OPM Flow scripts. The models are trained with the Keras library in Python scripts, stored in a format readable for the OPM Flow-NN framework and subsequently deployed within OPM Flow. To translate NNs from Keras to C++ functions, we  build on the existing tools Kerasify \citep{RobertW.RosePaulMaevskikh:Kerasify2017} and Keras2cpp \citep{RobertW.RoseGeorgyPerevozchikov:Keras2cpp2019}. Both are tools that convert deep learning models created within the Keras library to C++ code. The latter has, e.g., been used in the context of reservoir simulations for a \emph{Volume of Fluid-Machine Learning} (VOF-ML) scheme, which integrates ML flux functions into a finite volume solver \citep{Despres:MachineLearningDesign2020}. We modify and extend the scope of both tools to facilitate the deployment in OPM Flow and ensure that models are fully integrated into OPM Flow's automatic differentiation framework. When the user initializes and loads a stored NN inside an OPM Flow script, an automated deployment process handles all the translation. This process works by operating a series of steps handling model interpretation, layer conversion, optimization, and code generation steps to adapt the Keras model to a native OPM Flow function. Concretely, it consists of the following steps:
	\begin{itemize}
 		\item Model Parsing: First, the Keras model file generated after the training of the model is read and interpreted. The structure and layers of the model are analyzed.
   		\item Layer Conversion: Each layer of the Keras model is then converted into its equivalent C++ representation.
   		\item Optimization: Various optimizations are applied to improve the efficiency of the converted C++ code. This may involve combining certain layers, reducing unnecessary computations, or simplifying the network architecture.
   		\item Code Generation: Based on the parsed model and optimized settings, C++ code is generated using templating. This code is tailored to run efficiently on the target platform, taking advantage of available hardware acceleration, parallelization, automatic differentiation, or other optimizations specific to the platform.
   		\item Integration and Customization: The generated C++ code can then be incorporated into OPM Flow scripts by accessing it as a native function.
	\end{itemize}
	All but the last step are handled by the OPM Flow-NN framework and invisible to the user. This ensures easy handling. As shown in Listing \ref{listing:c++}, only a handful of lines are needed to access and deploy a NN in OPM Flow. 

	\begin{lstlisting}[caption={Example code on how to incorporate a stored NN into an OPM Flow script. }, language=C++, label=listing:c++]
KerasModel<Value> model; 	% Declare the NN model
model.LoadModel(path);		% Load a saved model
Tensor<Value> in{3};		% Declare the input tensor
in.data = {{p, re, V}};		% Initialize the input tensor
Tensor<Value> out;		% Declare the output tensor
model.Apply(&in, &out);		% Run the model
	\end{lstlisting}

	\section{Machine-learned near-well model}\label{section:near well model}
		
    The OPM Flow-NN framework is applied to a novel machine-learned near-well model. Similar to Peaceman-type and upscaling near-well models, a well index \(WI\) is prescribed that relates well flow rate, well pressure, and well block pressure (see Equation \eqref{equation:peaceman-original-model:linear-relation}). As in existing upscaling approaches, the well index is not derived (semi-)analytically, but upscaled from fine-scale simulations of near-well regions. Therefore, near-well effects are explicitly and accurately resolved. To allow for upscaling of complex transient and multiphase processes, an ensemble of fine-scale simulations with different simulation parameters is run. This ensemble approach enables capturing of various different flow regimes/geometries/etc. Afterward, a NN \(\mathcal{N}\) is employed as a means to store the upscaled ensemble data. More explicitly, the NN \(\mathcal{N}\) learns the map
    \begin{equation}\label{equation:nn map}
        \mathcal{D} \subset \R^n \to \R, \bm{x} \mapsto \widetilde{WI},
    \end{equation}
    where
    \begin{equation}\label{equation:data-driven WI}
    	\widetilde{WI} = \frac{q}{p_{well} - p_i}
    \end{equation} denotes the flow-based upscaled well index (in distinction to the analytic well index \eqref{equation:peaceman-original-model:WI}), calculated from the simulation values \(q,\) \(p_{well},\) and \(p_i.\) The vector \(\bm{x} = (x_1,\dots,x_n)\) is a set of simulation values that serve as inputs to the NN. These inputs can be both static and dynamical and shall be chosen to reflect the complexity of the flow regime and reservoir geometry - examples are given in section "\nameref{section:simulations}". In comparison to existing upscaling approaches \citep{Ding:ScalingupVicinityWells1995,Wolfsteiner:NearWellRadialUpscaling2002,Nakashima:AccurateRepresentationNearwell2010}, the NN approach conceptually allows capturing vastly more complexity, being limited only by the available ensemble data and the expressibility of the NN. If the ensemble data is sufficiently dense and the network sufficiently expressive, the accuracy of LGR \citep{Karimi-Fard:AccurateResolutionNearWell2012} is retained, but the computational load is shifted from during the simulation to before the simulation. Both the computational and conceptual difficulty lies then in constructing and running the ensemble simulations and training the network. If sufficiently many wells are represented by the same flow regimes and full-field simulations are run many times with the same NN near-well model (e.g., many simulations of the same field in optimization problems), a significant speedup in comparison to LGR may be achieved.

\begin{figure}[ht!]
        \centering
        \includeinkscape[inkscapelatex=false,width=\linewidth]{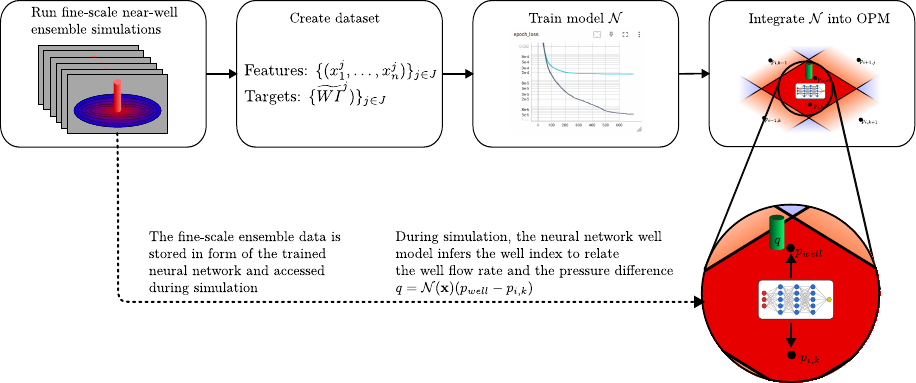}
        \caption[Workflow]{The recurring workflow to fit the machine-learned near-well model to previously unknown flow regimes and geometries.}
        \label{figure:workflow}
    \end{figure}

	A strength of this approach is that the workflow remains identical when incorporating different physical phenomena and reservoir geometries. They are simply included in the ensemble simulations and the datasets, and a new NN is trained. This makes the approach highly flexible and gives the capability to model complex relations between different factors that can influence the injectivity/productivity of a well. We showcase examples in section "\nameref{section:simulations}". The full workflow is illustrated in \Cref{figure:workflow}.

    We stress that, conceptually, the method is independent of the specific formulation of flow in porous media and can be incorporated in any reservoir simulator that prescribes a well index.

    \subsection{Governing equations - Black-oil formulation and well model}\label{subsection:black-oil-model}
		For demonstration purposes, we use the black-oil formulation \citep{Lake:FundamentalsEnhancedOil2014, Chen:ComputationalMethodsMultiphase2006} employed by OPM Flow. This formulation is an industry-standard for three-phase flow problems and implemented in various academic and commercial simulators. In the black-oil formulation, flow and transport of three phases and three pseudo-components, gas, oil and water, are modeled. The numerical examples in this work focus on the modeling of CO\textsubscript{2} storage. In these cases, the oil phase is prescribed the properties of brine, the gaseous phase the properties of CO\textsubscript{2}, and the water phase is disregarded. Gas can dissolve into the brine ("oil") phase and brine ("oil") can vaporize into the gas phase. The mass conservation equations in terms of the phase saturations \(S_\alpha\) and pressures \(p_\alpha, \, \alpha \in \{b, g\}\) (brine, gas/CO\textsubscript{2}) read as
        \begin{subequations}\label{equation:mass conservation}
        	\begin{align}
            	\frac{\partial}{\partial t} \left[\phi (\rho_b S_b + r_{bg} \rho_g S_g)\right] + \nabla \cdot (\rho_b \bm{v}_b + r_{bg} \rho_g \bm{v}_g) &= q_g, \\
            	\frac{\partial}{\partial t} \left[\phi (\rho_g S_g + r_{gb} \rho_b S_b)\right] + \nabla \cdot (\rho_g \bm{v}_g + r_{gb} \rho_b \bm{v}_b) &= q_b,
        	\end{align}
        \end{subequations}
        where \(\phi\) is porosity, \(\rho_\alpha\) is phase density, \(q_\alpha\) is the phase mass source/sink term, \(r_{bg}\) is the ratio of dissolved gas to brine and \(r_{gb}\) is the ratio of vaporized brine to gas. The phase fluxes are given by Darcy's law
        \begin{equation}\label{equation:Darcy's law}
            \bm{v}_\alpha = -\lambda_\alpha \bm{k} (\nabla p_\alpha - \rho_\alpha \bm{g}),
        \end{equation}
        where \(\lambda_\alpha = \frac{k_{r,\alpha}}{\mu_\alpha}\) is phase mobility, \(\bm{k}\) is the permeability tensor, \(\bm{g}\) is gravity, \(k_{r,\alpha}\) is relative permeability, and \(\mu_\alpha\) is phase viscosity. The equations are closed by the saturation constraint and the capillary pressure
        \begin{equation}\label{equation:closure terms}
            S_b + S_g = 1, \quad p_{c,bg} = p_b - p_g.
        \end{equation}
        Any production and injection wells can either be modeled explicitly by refining the grid in proximity, or implicitly by relating wellbore pressure, grid pressure, and flow by a well index \( WI.\) The exact derivation of \( WI\) depends on the employed near-well model. In the latter case, three additional variables per well are introduced to determine in-/outflow from wells \citep{Holmes:EnhancementsStronglyCoupled1983, Rasmussen:OpenPorousMedia2021}: the weighted total flow rate \(Q_t = \sum_{\alpha \in \{b, g\}} g_\alpha Q_\alpha,\) where \(Q_\alpha\) is component surface flow rate and \(g_\alpha\) is a weighting factor, the weighted gas fraction \(F_g = \frac{g_g Q_g}{Q_t},\) as well as the well's bottom hole pressure \(p_{bhp}.\) For each well \(w\) with a given set of connections \(C(w),\) the flow rate of phase \(\alpha\) at connection \(i\) is then modeled as
        \begin{equation}
            q_{\alpha, i}^r = WI_{\alpha,i} \left[p_i - (p_{bhp,w} + h_{w,i})\right],
        \end{equation}
        where \(p_i\) is the pressure value of the grid cell containing connection \(i\) and \(h_{w,i}\) is the difference between the bottom hole pressure and at the connection itself. Usually in OPM Flow, \( WI_{\alpha,i}\) is evaluated by a multiphase Peaceman model. In our new model, \( WI_{\alpha,i}\) is evaluated by a NN. The system is closed with mass balance equations for each component
        \begin{equation}
            \frac{A_{\alpha, w} - A_{\alpha, w}^0}{\Delta t} + Q_\alpha - \sum_{i \in C(w)} q_{\alpha, i} = 0,
        \end{equation}
        where \(A_{\alpha, w}\) denotes the component mass inside the well itself and \(q_{\alpha, i}\) is the component flow rate under surface conditions calculated from \(q_{\alpha, i}^r.\) For wells with pressure control, \(p_{bhp,w}\) is fixed, while for rate-controlled wells, \(Q_\alpha\) is fixed.

    \subsection{Fine-scale ensemble simulations and data generation}
    	To train the NN well model, three datasets \(D_{train},\) \(D_{val},\) and \(D_{test}\) for training, validation, and testing are constructed from the ensemble runs. Each consists of the model's inputs and the data-driven well index as a target
    	\begin{subequations}\label{equation:datasets}
        	\begin{align}
            	D_\alpha &= \{(x_1^j, \dots, x_n^j; \widetilde{WI}^j)\}_{j \in J_\alpha}, \\       
            	\widetilde{WI}^j &= \frac{q^j}{p_{well}^j - p_{i}^j}.
        	\end{align}
        \end{subequations}
        For the dataset to be an accurate representation of near-well phenomena, the governing \Cref{equation:mass conservation,equation:Darcy's law,equation:closure terms} are solved on a fine-scale grid in proximity of the wells without any additional well model. An ensemble of simulations with varying simulation parameters is run to capture varying flow regimes, transient effects, and differing reservoir geometries.

        While a NN can realistically learn the map \eqref{equation:nn map} for values inside the training data ranges for \(x_1,\dots,x_n,\widetilde{WI},\) valuable inference outside these ranges cannot be expected. Thus, the training data needs to cover all regimes and geometries the network is supposed to be employed on. For accurate interpolation between data points, the training data must also cover these ranges reasonably densely.

        These reasonings give general guidelines for how to construct the near-well ensemble simulations as a base for the datasets, but some thought may have to be put into specific cases. Examples are shown in section "\nameref{section:simulations}". For now, we denote the ensemble as follows. Let \(\nu_1 , \dots , \nu_m\) be the set of simulation input parameters that vary during different ensemble runs. The set of all ensemble members is then given by
        \begin{equation}\label{equation:ensemble}
            E = \{(\nu_1^i , \dots , \nu_m^i)\}_{i \in I},
        \end{equation}
        where \(|E| = |I|\) is the number of ensemble members.

\begin{figure}[h!]
            \centering
            \begin{minipage}[c]{0.67\textwidth}
            	\includeinkscape[inkscapelatex=false, width=\linewidth]{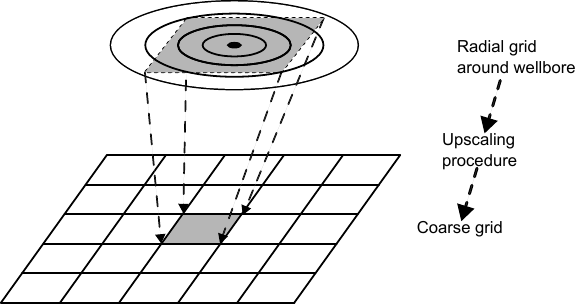}
            \end{minipage}\hfill
            \begin{minipage}[c]{0.3\textwidth}
            	\caption[Upscaling.]{Schematic of the upscaling procedure from the radial fine-scale simulations to a coarse-scale Cartesian cell.}
            	\label{figure:upscaling procedure}
            \end{minipage}
        \end{figure}

        To infer the correct values at simulation time, both the target well index \(\widetilde{WI}\) and the input values \(x_1,\dots,x_n\) are upscaled from the fine-scale, radial ensemble simulations to the coarse grid on which the model shall be employed on. For this work, we suppose that the coarse grid consists of square/square rectangular cuboid cells in the near-well region and that the fine-scale ensemble simulations are run on a radial grid (see \Cref{figure:upscaling procedure}). The upscaling procedures for relevant inputs for the examples in section "\nameref{section:simulations}" are then as follows:
        \begin{itemize}
            \item \(S_\alpha\) - phase saturation is arithmetically averaged over all fine-scale cells corresponding to a coarse-scale block. The fine-scale values are weighted by the overlapping fractions of the cells, i.e.,
            \begin{equation}\label{equation:saturation upscaling}
                S^{coarse}_{\alpha,i} = \sum\limits_{j \in J} f(ov)_{i,j} S^{fine}_{\alpha,j},
            \end{equation}
            where \(i\) is the index of the coarse cell, \(J\) is the set of all fine-scale cells, and \(f(ov)_{i,j} = \frac{Area(V_{j}^{fine} \cap V_{i}^{coarse})}{Area(V_{j})}\) is the fraction of fine cell \(V_j\) that overlaps with coarse cell \(V_i.\)
            \item \(p\) - as observed by Peaceman \citep{Peaceman:InterpretationWellBlockPressures1978}, the single pressure value corresponding to a coarse cell during simulation differs from the average pressure within the cell. For well blocks, where this phenomenon is particularly expressed, the cell pressure value is instead given by the pressure at equivalent well radius \(r_{e}.\) Hence, the upscaled value is taken at the radial cell whose radius is closest to \(r_e\):
                \begin{equation}\label{equation:pressure upscaling}
                    p^{coarse}_{\alpha,i} = p^{fine}_{\alpha,j} \text{ s.t. } j = \argmin_k \|r_e - r_k\|, \text{ where } r_k \text{ is the radius of fine-scale cell }k.
                \end{equation}
            \item \(V_{tot},\dots\) - global variables such as total injected volume are left unchanged.
        \end{itemize}
		For more complicated geometries, we refer to \citep{Ding:ScalingupVicinityWells1995,Wolfsteiner:NearWellRadialUpscaling2002,Nakashima:AccurateRepresentationNearwell2010} for upscaling procedures.

        The domain size for the ensemble simulations is determined by the desired size of the coarse-size well cells that the model shall be employed on. To relate the well block pressure with the wellbore pressure, the ensemble simulations need to encompass all pressure isobars up to at least the one with pressure value equal to the coarse block pressure. The radius of this isobar depends on the coarse-scale grid size. Assuming homogeneity inside the near-well region, it is given by the equivalent well radius \(r_e\) and can, e.g., be found by the techniques described in \citep{Peaceman:InterpretationWellBlockPressures1978,Peaceman:InterpretationWellBlockPressures1983, Peaceman:InterpretationWellblockPressures1990,Aavatsmark:WellIndexReservoir2003}. When the model inputs contain simulation data that is upscaled by averaging over the entire region enclosed by the coarse grid block, the fine-scale domain needs to have according domain size. If the model inputs contain information from neighboring coarse grid blocks (e.g., example \nameref{subsection:CO2 3D example}), the according domain needs to be included etc.

        Note that for each ensemble member, multiple data points are extracted, corresponding to different time steps and coarse-scale cells sizes, i.e., different radii at fine-scale level. Hence, in general \(|D_{train} \cup D_{val} \cup D_{test}| > |E|.\)

    \subsection{Model architecture and training}
\begin{figure}[h!]
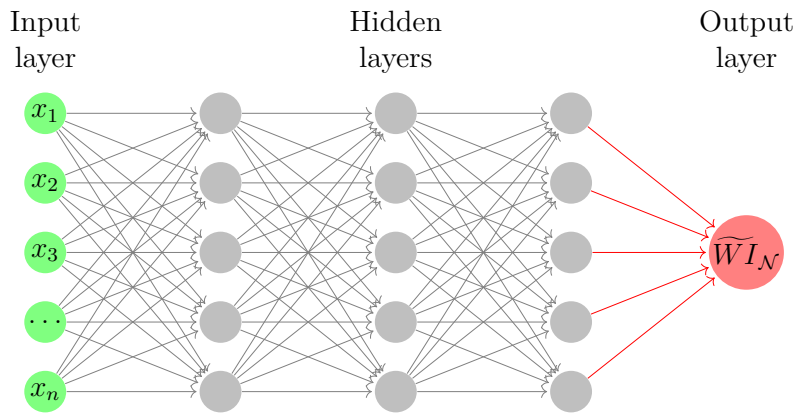

            \centering
            \begin{minipage}[c]{0.67\textwidth}        
                \resizebox{\linewidth}{!}{
            	\begin{neuralnetwork}
        			\newcommand{\inputtext}[2]{\ifcase#2 \relax \or \(x_1\) \or \(x_2\) \or \(x_3\) \or \(\dots\) \or \(x_n\) \fi}
        			\newcommand{\outputtext}[2]{\ifcase#2 \relax \or \(\widetilde{WI}_\mathcal{N}\) \fi}
        			\tikzset{linkfix/.style = {draw=gray} }
        			\tikzset{linkflex/.style = {draw=red} }
        	       	\inputlayer[count=5, title=Input\\layer, bias=false, text=\inputtext]
        			\layer[count=5, bias=false, nodeclass={neuron}]
        	       	\linklayers[style={linkfix}]
        			\layer[count=5, bias=false, nodeclass={neuron},title=Hidden\\layers]
           		   	\linklayers[style={linkfix}]
        			\layer[count=5, bias=false, nodeclass={neuron}]
                	\linklayers[style={linkfix}]
               	 	\outputlayer[count=1, title=Output\\layer, text=\outputtext]
                	\linklayers[style={linkflex}]
            	\end{neuralnetwork}
            	}
            \end{minipage}\hfill
            \begin{minipage}[c]{0.3\textwidth}
            	\caption[Neural network.]{Scheme of the FCNN model. Inputs can be chosen freely to accustom for the complexity of the near-well flow regimes and the geometry.}
            	\label{figure:neural network}
            \end{minipage}
        \end{figure}

        A fully connected NN (FCNN) (see \Cref{figure:neural network}) is employed to learn the map \eqref{equation:nn map}. Before training, a hyperparameter search is performed to find the best architecture, i.e., depth, hidden size, and activation function of the network.  The loss function is the mean-squared error between the predicted \(\widetilde{WI}_\mathcal{N}^j\) and the target well index \(\widetilde{WI}^j\)
        \begin{equation}
            \mathcal{L} = \frac{1}{K} \sum_{k = 1}^K |\widetilde{WI}_\mathcal{N}^k - \widetilde{WI}^k|^2,
        \end{equation}
        where \(K\) is the batch size (i.e., by abuse of notation, \(\{1,\dots,K\} \subset D_\alpha\)). Note, that all inputs and outputs in \(D_{train} \cup D_{val} \cup D_{test}\) are scaled to the interval \([-1,1]\) with min-max normalization. The network is optimized with Adam \citep{Kingma:AdamAMethodforStochasticOptimization} until the validation loss stays constant for sufficiently many epochs. We emphasize that the workflow is not limited to any specific NN architecture. For complex flow regimes and/or geometries, it may be beneficial to employ more sophisticated architectures than FCNN.

	\section{Numerical simulations}\label{section:simulations}
		
    To showcase the flexibility of the method, we train three models in different flow regimes and on geometries that vary in complexity. The second example (subsection "\nameref{subsection:CO2 2D example}") highlights the capability to include expert knowledge into the model, while the third example (subsection "\nameref{subsection:CO2 3D example}") highlights the ability of the NN approach to accurately model complex relations without the need for mathematical analysis. The models' accuracy is tested on reservoir-scale simulations against a fine-scale benchmark with LGR in well vicinity. The results are compared to the accuracy of multiphase Peaceman correction. OPM Flow features solely the multiphase Peaceman method, hence providing a model with improved accuracy enhances OPM Flow's modeling capacities significantly. Note that all examples are run with rate control for the wells and a fixed injection rate. The examples focus on representing buoyancy and phase-flow effects on the well injectivity with a NN, hence dissolution of the CO\textsubscript{2} is not considered. The parameters for the CO\textsubscript{2} examples are taken from simulations of real-world CO\textsubscript{2} injection in the Sleipner field in the North Sea \citep{Singh:ReservoirModelingCO22010}, but injection rate and anisotropy ratio are increased to emphasize the near-well effects.

    \subsection{Single-phase H\texorpdfstring{\textsubscript{2}}{2}O injection in homogeneous two-dimensional reservoir}\label{subsection:H2O example}
        To verify the validity of the framework and to affirm that the model can learn a basic inflow-performance relation, we consider the injection of H\textsubscript{2}O into a homogeneous, two-dimensional, water filled reservoir (the domain is assigned a homogeneous height over one vertical cell). Under these conditions, the governing equations \cref{equation:mass conservation,equation:Darcy's law,equation:closure terms} reduce to single-phase flow and the Peaceman model will give the exact relationship between well flow, wellbore pressure, and cell pressure on a repeated five-spot pattern on a Cartesian grid \citep{Peaceman:InterpretationWellBlockPressures1978}. For the single-well injection on a Cartesian grid that we run, the Peaceman model gives an accurate solution on all resolutions.

        The Peaceman formula for the well index \eqref{equation:peaceman-original-model:WI} gives a natural choice of model inputs and ensemble parameters. The inflow performance relationship is dependent on \(k,\) reservoir height \(h,\) \(r_e,\) \(\rho(p),\) and \(\mu(p).\) We take the former three directly as inputs for the NN and use cell pressure \(p\) as a notion of density and viscosity. The model's input vector is then \(\left(x_1, x_2, x_3, x_4\right) =  \left(p, k, h, r_e\right),\) all values taken from the grid cell enclosing the well. The ensemble is constructed as \(E = \{\left(p_{init}^i, k^i, h^i\right)\}_{i \in I}\) s.t. each of the inputs is densely represented by multiple data points inside some range.

\begin{figure}[p!]
            \centering
            \begin{minipage}[c]{0.67\textwidth}
                \centering
                \includeinkscape[inkscapelatex=false, width=\linewidth]{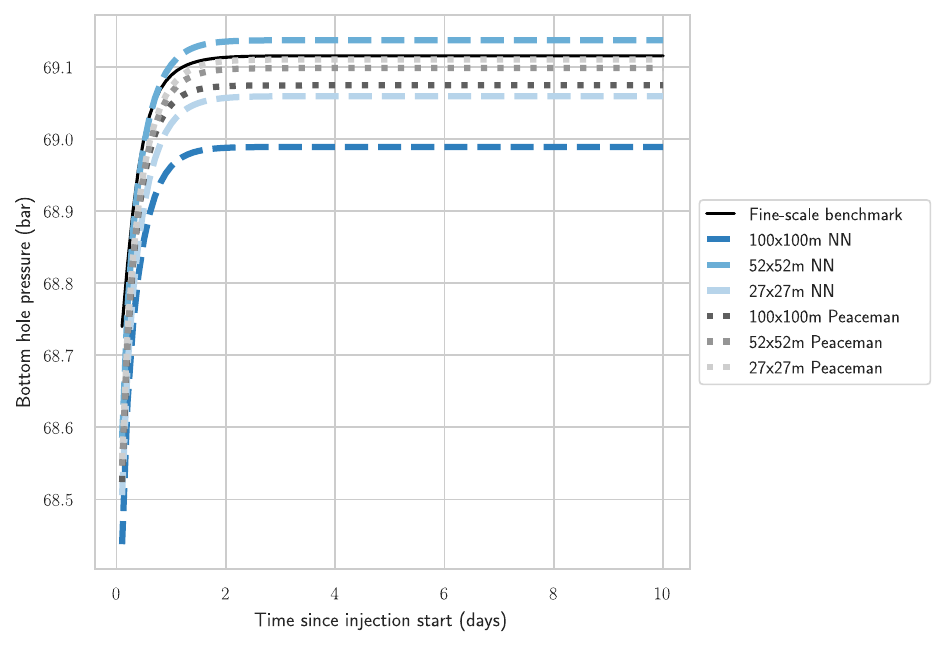}
            \end{minipage}\hfill
            \begin{minipage}[c]{0.3\textwidth}
                \subcaption{Simulation values: \\ \(k = \qty{2d-13}{\meter\squared}, p_{init} = \qty{65}{\bar}, h = \qty{7.5}{\meter}\)}
                \label{figure:H2O bhp_1}
            \end{minipage}
            \bigskip
            \begin{minipage}[c]{0.67\textwidth}
                \centering
                \includeinkscape[inkscapelatex=false, width=\linewidth]{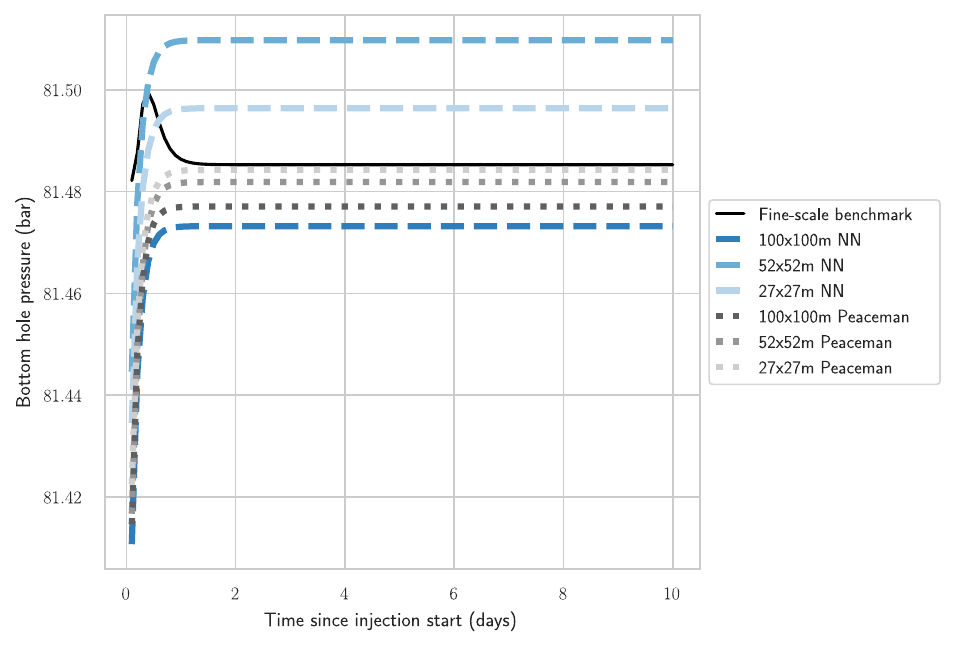}
            \end{minipage}\hfill
            \begin{minipage}[c]{0.3\textwidth}
                \subcaption{Simulation values: \\ \(k = \qty{5d-13}{\meter\squared}, p_{init} = \qty{80}{\bar}, h = \qty{15}{\meter}.\) Note that \(p_{bhp}\) for all simulations stays inside a range of only \(\qty{0.15}{\bar},\) i.e., all simulations are nearly exact despite \(p_{bhp}\) decreasing in the fine-scale benchmark and increasing in coarser simulations.}
                \label{figure:H2O bhp_2}
            \end{minipage}
			\bigskip
            \begin{minipage}[c]{0.67\textwidth}
                \centering
                \includeinkscape[inkscapelatex=false, width=\linewidth]{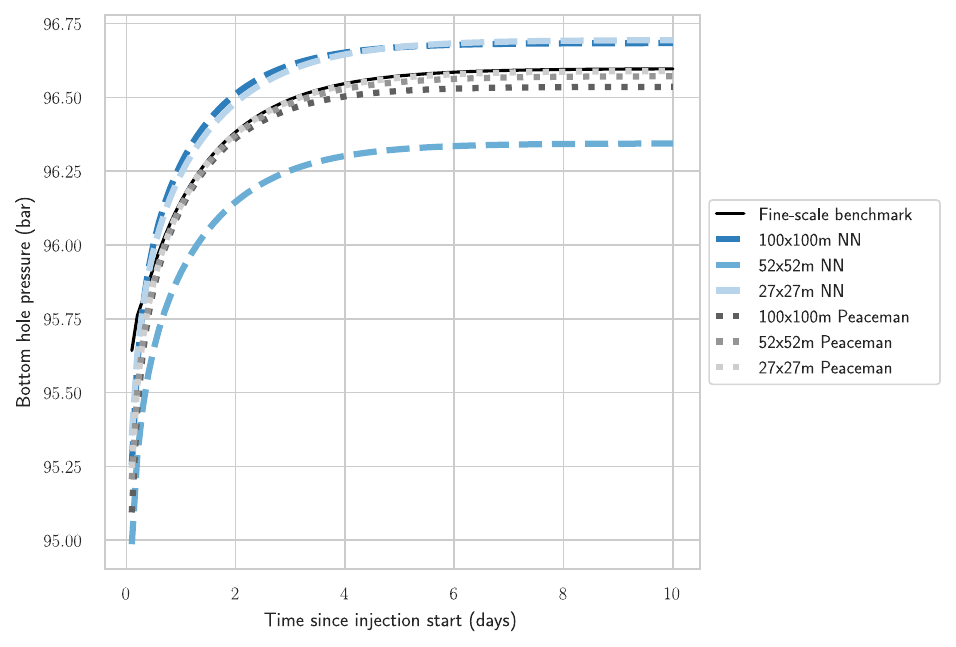}
            \end{minipage}\hfill
            \begin{minipage}[c]{0.3\textwidth}
                \subcaption{Simulation values: \\ \(k = \qty{5d-14}{\meter\squared}, p_{init} = \qty{90}{\bar}, h = \qty{20}{\meter}\)}
                \label{figure:H2O bhp_3}
            \end{minipage}
            \caption[H\textsubscript{2}O \(p_{bhp}\)]{H\textsubscript{2}O example: \(p_{bhp}\) inferred by the ML model and the Peaceman model on different mesh sizes for three different simulations. Peaceman gives the exact solution independent of resolution, which the ML model is able to closely match for all simulations and mesh sizes.}
            \label{figure:H2O bhp}
        \end{figure}

        The trained model is compared to Peaceman on three reservoir-scale simulation, which differ in initial pressure, reservoir permeability, and domain height. Input parameters for the ensemble and full-scale simulations are shown in \cref{table:H2O parameters} and \cref{figure:H2O bhp}. Note that all simulation parameters the NN was not explicitly trained on, such as wellbore radius and temperature, stay unchanged between ensemble and final simulation. The network learns these values implicitly. In \cref{figure:H2O bhp}, we compare the bottom hole pressure \(p_{bhp}\) over time inferred by the NN well model and multiphase Peaceman correction. Peaceman infers bottom hole pressure values close to the fine-scale benchmark values independent of the grid size, confirming the theoretical result that it gives an accurate representation of the analytical solution for radial single-phase flow in homogeneous rock  \citep{Peaceman:InterpretationWellBlockPressures1978}. The ML model is accurate up to maximal errors of \(\qty{0.4}{\bar},\) \(\qty{0.1}{\bar},\) and \(\qty{0.7}{\bar}\) and average errors of \(\qty{0.2}{\bar},\) \(\qty{0.1}{\bar},\) and \(\qty{0.3}{\bar}\) for the three simulations, respectively, over the time of injection. This assures that the ML model is able to learn the basic inflow-performance relation \eqref{equation:peaceman-original-model:WI} and can provide accurate results during simulation in OPM Flow.

        \begin{table}[ht!]
            \centering
            \begin{tabularx}{.9\textwidth}{|>{\raggedright\arraybackslash}X p{1.5cm} >{\raggedright\arraybackslash}X >{\raggedright\arraybackslash}X|}
                \hline
                parameter &  symbol & ensemble values & full simulation values \\
                \hline
                horizontal domain size & \(-\)  & \(\qtyproduct{100}{\meter}\) &  \(\qtyproduct{1100 x 1100}{\meter}\) \\
                vertical domain size & \(h\)  & \(\qtyrange[range-units = single, range-phrase =-]{2}{20}{\meter}\) &  \(\qtyrange[range-units = single, range-phrase =-]{7.5}{20}{\meter}\) \\
                cell sizes & \(-\)  & log. grid, 50 cells of varying sizes in \(x\)-direction & \(\qtyproduct{100 x 100}{\meter}\) (coarse) to \(\qtyproduct{27 x 27}{\meter}\) (fine) \\
                porosity    & \(\phi\)  & \(\qty{0.35}{}\) & \(\qty{0.35}{}\) \\
                permeability    & \(k\)  & \(\qtyrange[range-units = single, range-phrase =-]{1d-14}{1d-12}{\meter\squared}\) & \(\qtyrange[range-units = single, range-phrase =-]{5d-14}{5d-13}{\meter\squared}\) \\
                injection rate    & \(Q_t\)  & \(\qty{60}{\meter\cubed\per\day}\) & \(\qty{60}{\meter\cubed\per\day}\) \\
                initial reservoir pressure    & \(p_{init}\)  & \(\qtyrange[range-units = single, range-phrase =-]{50}{150}{\bar}\) & \(\qtyrange[range-units = single, range-phrase =-]{65}{90}{\bar}\) \\
                reservoir temperature & \(T\)  & \(\qty{40}{\degreeCelsius}\) & \(\qty{40}{\degreeCelsius}\) \\
                wellbore radius & \(r_w\)  & \(\qty{0.25}{\meter}\) & \(\qty{0.25}{\meter}\) \\
                \hline
            \end{tabularx}
           	\captionsetup{width=.9\textwidth}
            \caption[H\textsubscript{2}O parameters]{Simulation parameters for the H\textsubscript{2}O example. Note that both the ensemble domain and the full domain have large pore-volume cells at the horizontal boundaries.  For the fine-scale benchmark, cell sizes are \(\qtyproduct{4.5 x 4.5}{\meter}\) in well vicinity and \(\qtyproduct{18 x 18}{\meter}\) further out.}
            \label{table:H2O parameters}
        \end{table}

    \subsection{Two-phase CO\texorpdfstring{\textsubscript{2}}{2} injection in homogeneous two-dimensional reservoir}\label{subsection:CO2 2D example}
        In the second example, we consider the injection of CO\textsubscript{2} into a homogeneous, two-dimensional, water-saturated domain. For multiphase flow, the near-well solution is, in general, not steady-state. Therefore, in addition to the inputs from the previous example, the NN requires some notion of the reservoir state as an input. We solve this by passing the total injected volume of CO\textsubscript{2}, denoted by \(V_{tot},\) to the network. As the injection rate is constant in both the ensemble and full reservoir simulations, \(V_{tot}\) in combination with grid block pressure \(p_i\) suffices to characterize the system's state to good accuracy. As noted, the network's inputs and outputs are scaled to \([-1,1].\) Furthermore, some of the feature and output data, such as \(\widetilde{WI},\) spreads over an exponential range. When the network infers \(WI\) close to the minimum of the range, small perturbations in the inputs might lead to network outputs just below \(-1,\) which correspond to negative and thus nonphysical \(\widetilde{WI}\) when scaled back. To prevent such behavior and improve accuracy, the model is trained to infer \(\log_{10} \left(\widetilde{WI}\right)\) instead for both this and the next example.

        \subsubsection{Including expert knowledge in the model}
        	Making use of the flexibility of the NN approach, expert knowledge can be passed directly to the model. This can, e.g., be the analytical near-well models discussed in sections "\nameref{section:introduction}" and "\nameref{section:simulations}". By including expert knowledge directly, it does not have to be learned by the NN, which enables faster training and better accuracy. We showcase this by passing the geometric portion of the Peaceman model \eqref{equation:peaceman-original-model:WI} as a network input. The model then effectively learns a correction term stemming from varying phase mobility in well proximity. Assuming radial flow and homogeneous permeability, the flow rate \(q_\alpha\) for phase \(\alpha\) at a single well connection can be computed in terms of the volume flux \(\vec{v}_\alpha\) by integrating around a circle \(\mathcal{C}\) of radius \(r\) around the connection
        	\begin{equation}
            	q_{\alpha} = \int_\mathcal{C} \rho_{\alpha}\vec{v}_{\alpha}\cdot\vec{n}ds = -2 \pi h r \rho_{\alpha} v_{\alpha}.
        	\end{equation}
        	Ignoring any potential vertical flow and integrating from \(r_w\) to \(r,\) one obtains
        	\begin{align}
            	\frac{2 \pi k h}{q_\alpha} \int_{r_{w}}^{r}\frac{k_{r,\alpha} \rho_{\alpha}}{\mu_{\alpha}} \frac{dp_\alpha}{dr} dr &= \int_{r_{w}}^{r} \frac{1}{r} dr. \\
            	\intertext{Approximating the interval on the left hand side with \(\frac{\tilde{k}_{r,\alpha} \tilde{\rho}_{\alpha}}{\tilde{\mu}_{\alpha}} \left[p_\alpha(r) - p_{\alpha,w}\right],\) where \(\tilde{k}_{r,\alpha},\) \(\tilde{\rho}_{\alpha},\) and \(\tilde{\mu}_{\alpha}\) are averaged values, one obtains}
            	\frac{2 \pi k h}{\ln{\frac{r}{r_w}}} \left(\frac{\tilde{k}_{r,\alpha} \tilde{\rho}_{\alpha}}{\tilde{\mu}_{\alpha}} + E\right) \left[p_\alpha(r) - p_{\alpha,w}\right] &= q_\alpha, \\
            	\intertext{with \(E\) being the error between the approximate and exact solution. Setting \(p_{\alpha,i} = p_\alpha(r_e)\) in equation \eqref{equation:peaceman-original-model:linear-relation}, it holds} \label{equation:correction factor}
            	WI = \frac{2 \pi k h}{\ln{\frac{r_e}{r_w}}} \left(\frac{\tilde{k}_{r,\alpha} \tilde{\rho}_{\alpha}}{\tilde{\mu}_{\alpha}} + E\right).
        	\end{align}
			By passing the logarithm of the geometrical part of the well index
			\begin{equation}\label{equation:expert knowledge}
				\log_{10} \left(\frac{2 \pi k h}{\ln{\frac{r_e}{r_w}}}\right)
			\end{equation}
			directly to the NN, it only has to learn the correction factor \(\log_{10} \left(\frac{\tilde{k}_{r,\alpha} \tilde{\rho}_{\alpha}}{\tilde{\mu}_{\alpha}} + E\right)\) to infer \(\log_{10} \left(\widetilde{WI}\right).\) Thus, \(k,\) \(h,\) and \(r_e\) are not passed directly to the NN, but implicitly as \(\frac{2 \pi k h}{\ln{\frac{r_e}{r_w}}}.\)  The full model input vector is \(\left(x_1, x_2, x_3\right) = \left(p, \log_{10} \left(\frac{2 \pi k h}{\ln{\frac{r_e}{r_w}}}\right), V_{tot}\right),\) with the local values again taken from the cell enclosing the well. The ensemble is constructed as \(E = \{\left(p_{init}^i, k^i\right)\}_{i \in I}.\)

\begin{figure}[p!]
            \centering
            \begin{minipage}[c]{0.67\textwidth}
                \centering
                \includeinkscape[inkscapelatex=false, width=\linewidth]{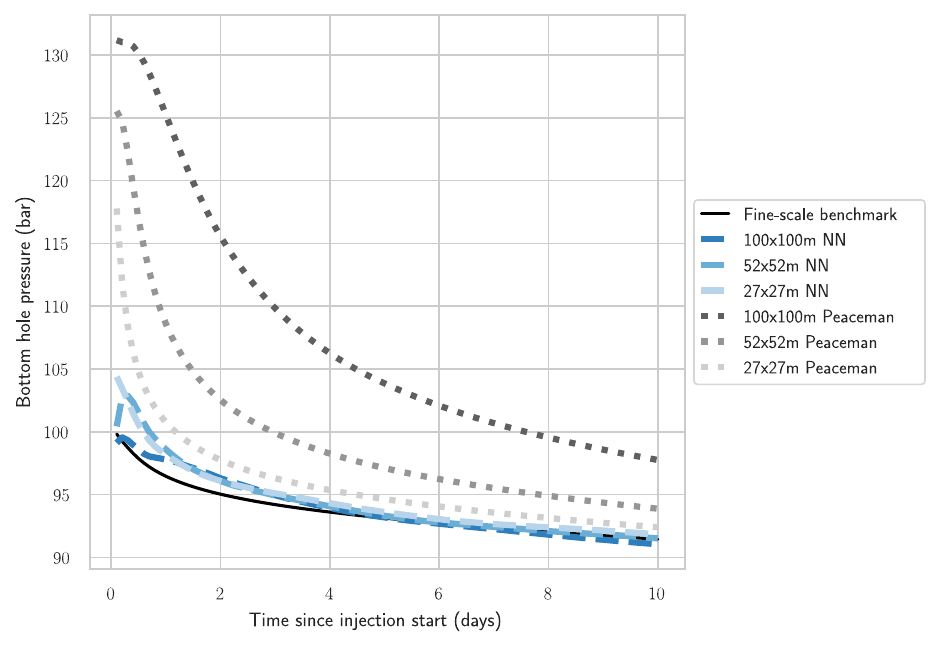}
            \end{minipage}\hfill
            \begin{minipage}[c]{0.3\textwidth}
                \subcaption{Simulation values: \\ \(k = \qty{1d-12}{\meter\squared}, p_{init} = \qty{65}{\bar}\)}
                \label{figure:CO2 2D bhp_1}
            \end{minipage}
            \bigskip
            \begin{minipage}[c]{0.67\textwidth}
                \centering
                \includeinkscape[inkscapelatex=false, width=\linewidth]{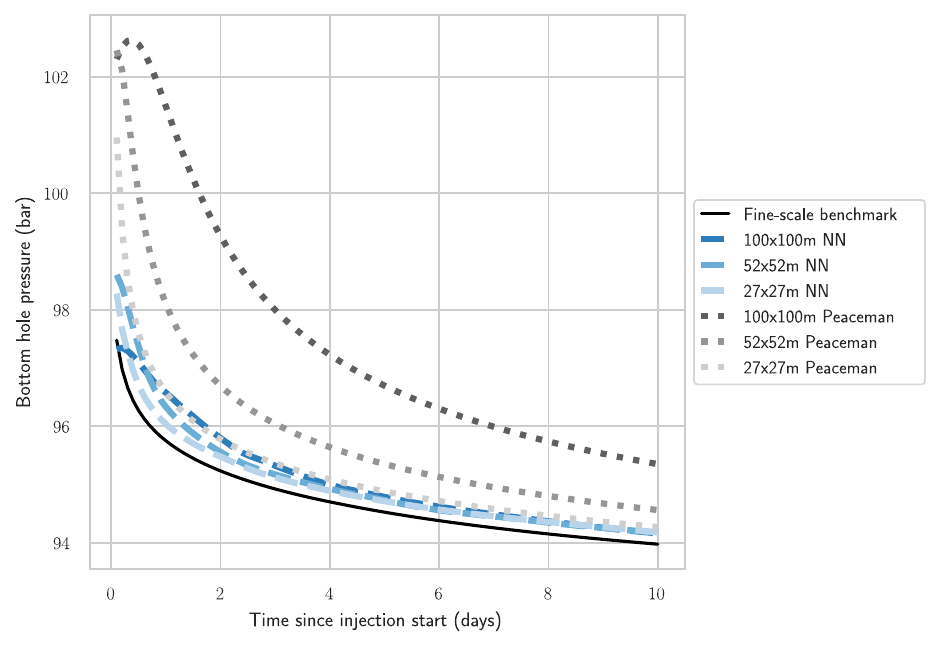}
            \end{minipage}\hfill
            \begin{minipage}[c]{0.3\textwidth}
                \subcaption{Simulation values: \\ \(k = \qty{5d-12}{\meter\squared}, p_{init} = \qty{90}{\bar}\)}
                \label{figure:CO2 2D bhp_2}
            \end{minipage}
			\bigskip
            \begin{minipage}[c]{0.67\textwidth}
                \centering
                \includeinkscape[inkscapelatex=false, width=\linewidth]{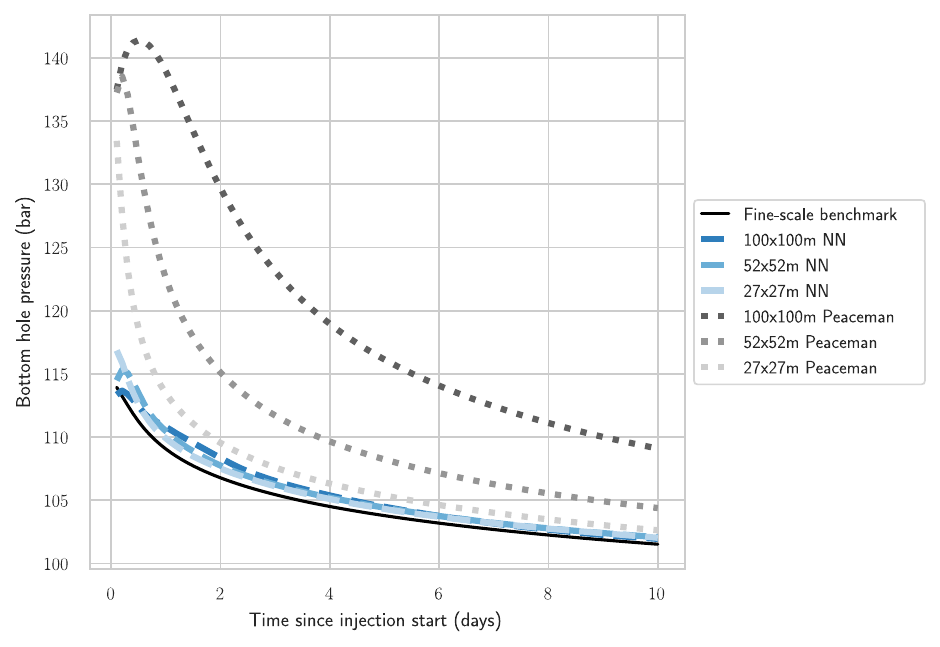}
            \end{minipage}\hfill
            \begin{minipage}[c]{0.3\textwidth}
                \subcaption{Simulation values: \\ \(k = \qty{9d-13}{\meter\squared}, p_{init} = \qty{80}{\bar}\)}
                \label{figure:CO2 2D bhp_3}
            \end{minipage}
            \caption[CO\textsubscript{2} \(p_{bhp}\)]{CO\textsubscript{2} 2D example: \(p_{bhp}\) inferred by the ML model and multiphase Peaceman on coarse meshes in comparison with a fine-scale benchmark solution. The ML model, integrated into a coarse-scale simulation, is able to accurately infer the fine-scale dynamics in all simulations.}
            \label{figure:CO2 2D bhp}
        \end{figure}

        Input parameters for the ensemble and full-scale simulations are shown in \cref{table:CO2 2D parameters}. We plot the bottom hole pressure over time for three reservoir-scale injection scenarios (\cref{figure:CO2 2D bhp}) with differing permeabilities and initial pressures, and compare to multiphase Peaceman on different grid scales. The ML model is accurate to the fine-scale benchmark up to maximal errors of \(\qty{4.6}{\bar},\) \(\qty{1.4}{\bar},\) and \(\qty{2.9}{\bar}\) and average errors of \(\qty{0.8}{\bar},\) \(\qty{0.4}{\bar},\) and \(\qty{0.8}{\bar}\) for the three simulations, respectively, over the time of injection. Notably, the model attains almost the same accuracy for all grid resolutions, demonstrating that it can infer fine-scale phenomena independent of grid size. In comparison, multiphase Peaceman gives inaccurate bottom hole pressure values for the reservoir-scale simulations. On the coarsest grid, the maximal errors are \(\qty{32.4}{\bar},\) \(\qty{6.3}{\bar},\) and \(\qty{30.8}{\bar}\) and the average errors \(\qty{13.8}{\bar},\) \(\qty{2.8}{\bar},\) and \(\qty{15.1}{\bar},\) all significantly higher than the ML model. Only for grid cells of size \(\qtyproduct{27 x 27 x 5}{\meter},\) multiphase Peaceman achieves a result comparable to the ML model run on the coarsest grid, i.e., Peaceman requires 14 times finer grid resolution for similar accuracy. Nevertheless, in the initial hours, the Peaceman model's error is still significantly higher than the ML model on all grid resolutions. This result highlights the capability of the ML model to accurately infer fine-scale transient effects from from a coarse-scale simulation at runtime.

        \begin{table}[ht!]
            \centering
            \begin{tabularx}{.9\textwidth}{|>{\raggedright\arraybackslash}X p{1.5cm} >{\raggedright\arraybackslash}X >{\raggedright\arraybackslash}X|}
                \hline
                parameter &  symbol & ensemble values & full simulation values \\
                \hline
                domain size & \(-\)  & \(\qtyproduct{100 x 5}{\meter}\) (radial) &  \(\qtyproduct{110 x 110 x 5}{\meter}\) \\
                cell sizes & \(-\)  & log. grid, 50 cells of varying sizes in \(x\)-direction &  \(\qtyproduct{100 x 100 x 5}{\meter}\) (coarse) to \(\qtyproduct{27 x 27 x 5}{\meter}\) (fine) \\
                porosity    & \(\phi\)  & \(\qty{0.35}{}\) & \(\qty{0.35}{}\) \\
                permeability    & \(k\)  & \(\qtyrange[range-units = single, range-phrase =-]{5d-13}{1d-11}{\meter\squared}\) & \(\qtyrange[range-units = single, range-phrase =-]{9d-13}{5d-12}{\meter\squared}\) \\
                injection rate    & \(Q_t\)  & \(\qty{9342.15}{\tonne\per\day}\) & \(\qty{9342.15}{\tonne\per\day}\) \\
                initial reservoir pressure    & \(p_{init}\)  & \(\qtyrange[range-units = single, range-phrase =-]{50}{120}{\bar}\) & \(\qtyrange[range-units = single, range-phrase =-]{65}{90}{\bar}\) \\
                reservoir temperature & \(T\)  & \(\qty{40}{\degreeCelsius}\) & \(\qty{40}{\degreeCelsius}\) \\
                wellbore radius & \(r_w\)  & \(\qty{0.25}{\meter}\) & \(\qty{0.25}{\meter}\) \\
                \hline
            \end{tabularx}
           	\captionsetup{width=.9\textwidth}
            \caption[CO\textsubscript{2} 2D parameters]{Simulation parameters for the 2D CO\textsubscript{2} example. Note that both the ensemble domain and the full domain have large pore-volume cells at the horizontal boundaries. For the fine-scale benchmark, cell sizes are \(\qtyproduct{4.5 x 4.5 x 5}{\meter}\) in well vicinity and \(\qtyproduct{18 x 18 x 5}{\meter}\) further out.}
            \label{table:CO2 2D parameters}
        \end{table}

    \subsection{Two-phase CO\texorpdfstring{\textsubscript{2}}{2} injection in vertically heterogeneous three-dimensional reservoir}\label{subsection:CO2 3D example}
\begin{figure}[hbt!]
            \centering
            \begin{subfigure}[t]{.45\textwidth}
                \centering
                \includegraphics[width=\linewidth]{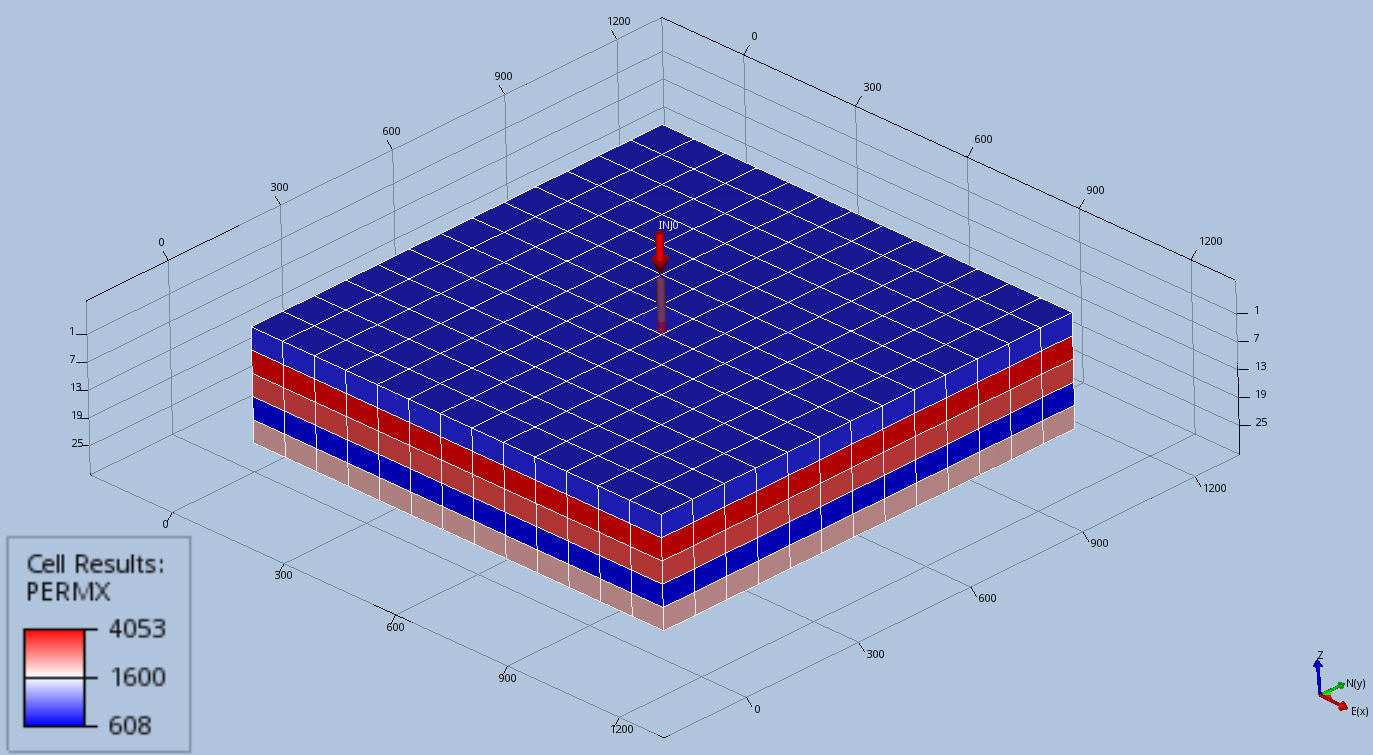}
                \subcaption{Geometry of the full-scale simulation. The vertical scale of the domain is exaggerated for illustration.}
                \label{figure:CO2 3D geometry}
            \end{subfigure}
            \begin{subfigure}[t]{.45\textwidth}
                \centering
                \includeinkscape[inkscapelatex=false, width=\linewidth]{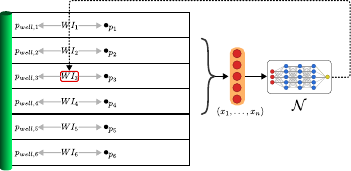}
                \subcaption[Local stencil]{The local stencil model takes into account coarse-cell simulation data not only from the well connection block itself, but also from neighboring cells.}
                \label{figure:local stencil}
            \end{subfigure}
            \caption{Exemplary geometry and ML model for the 3D CO\textsubscript{2} example with a layered reservoir.}
        \end{figure}

        Last, we consider a CO\textsubscript{2} injection in 3D with vertical heterogeneous layers, as shown in \cref{figure:CO2 3D geometry}. The anisotropy ratio of the domain is taken to be \(\frac{k_{xy}}{k_z} = 0.5,\) where \(k_{xy}\) is horizontal permeability and \(k_z\) is vertical permeability. Due to the density difference between CO\textsubscript{2} and water, vertical flow between layers will occur and influence the near-well pressure field. Therefore, the well index can depend on inflow from neighboring cells, in addition to the dependency on the well connection cell itself. To model this analytically would require significant work. The machine-learning approach makes modeling much more flexible, as we can simply input the relevant data into the NN and let it determine the correct mapping. The model inputs are therefore chosen to contain both data from the well connection grid block itself, but also from vertically connected blocks (see \cref{figure:local stencil}). This showcases the capability of the NN approach to detect patterns that are either imperceptible to normal observation and/or difficult to capture in mathematical terms.

		The full input vector to the network is
        \begin{equation}
            \left(x_1, x_2, \dots, x_{12}\right) = \left(p_u, p, p_l, S_{g,u}, S_g, S_{g,l}, k_{xy,u}, k_{xy}, k_{xy,l}, V_{tot}, r_e, \log_{10} \left(\frac{2 \pi k h}{\ln{\frac{r_e}{r_w}}}\right)\right),
        \end{equation}
        with the local values taken from the cell enclosing the well connection,  the cell above \((\cdot)_u,\) and the cell below \((\cdot)_l,\) respectively. While constructing the datasets, the data has to be padded s.t. the input vector is well-defined for upper and lower boundary cells. The pressure is padded with the values of the boundary cells, while permeability and saturation are padded with zero values. We highlight the capability of the NN approach to include a large number of different input data, which again would be challenging in an analytic approach.

        The ensemble simulations feature vertical heterogeneity with randomly chosen permeability for each layer and varying initial pressure. The ensemble is constructed as
        \begin{equation}
            E = \{(p_{init}^i, k_1^i, _2^i, \dots, k_L^i)\}_{i \in I},
        \end{equation}
        where \(L\) is the number of vertical layers of differing permeability.

\begin{figure}[p!]
            \centering
            \begin{minipage}[c]{0.67\textwidth}
                \centering
                \includeinkscape[inkscapelatex=false, width=\linewidth]{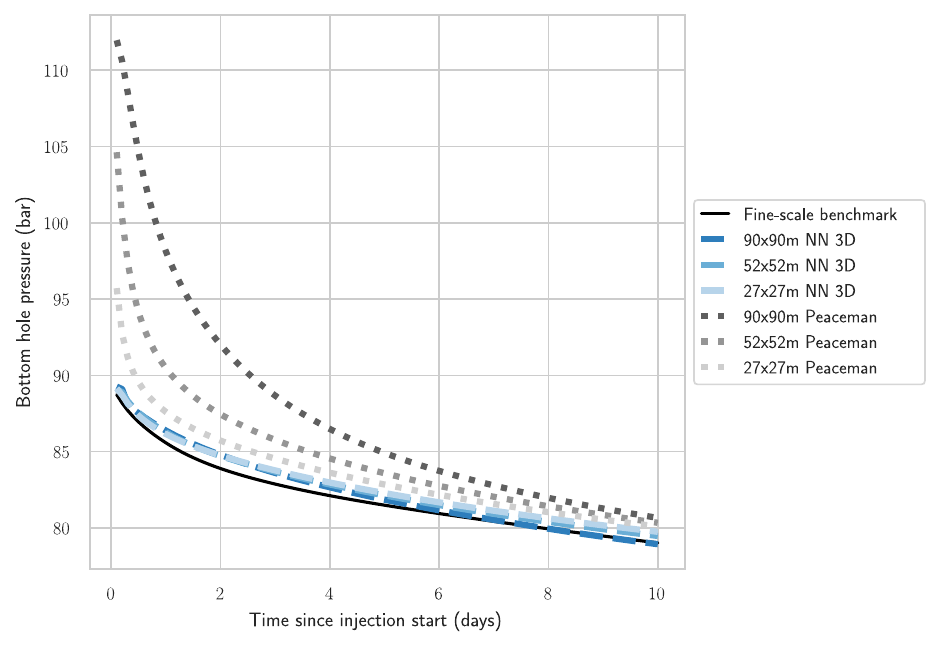}
            \end{minipage}\hfill
            \begin{minipage}[c]{0.3\textwidth}
                \subcaption[filler]{Simulation values:
                	\begin{align*}
                		k_{xy,1} &= \qty{7d-13}{\meter\squared} \\
                		k_{xy,2} &= \qty{4d-12}{\meter\squared} \\
                		k_{xy,3} &= \qty{3d-12}{\meter\squared} \\
                		k_{xy,4} &= \qty{6d-13}{\meter\squared} \\
                		k_{xy,5} &= \qty{2d-12}{\meter\squared} \\
                		p_{init} &= \qty{65}{\bar}
					\end{align*}           
                }
                \label{figure:CO2 3D bhp_1}
            \end{minipage}
            \bigskip
            \begin{minipage}[c]{0.67\textwidth}
                \centering
                \includeinkscape[inkscapelatex=false, width=\linewidth]{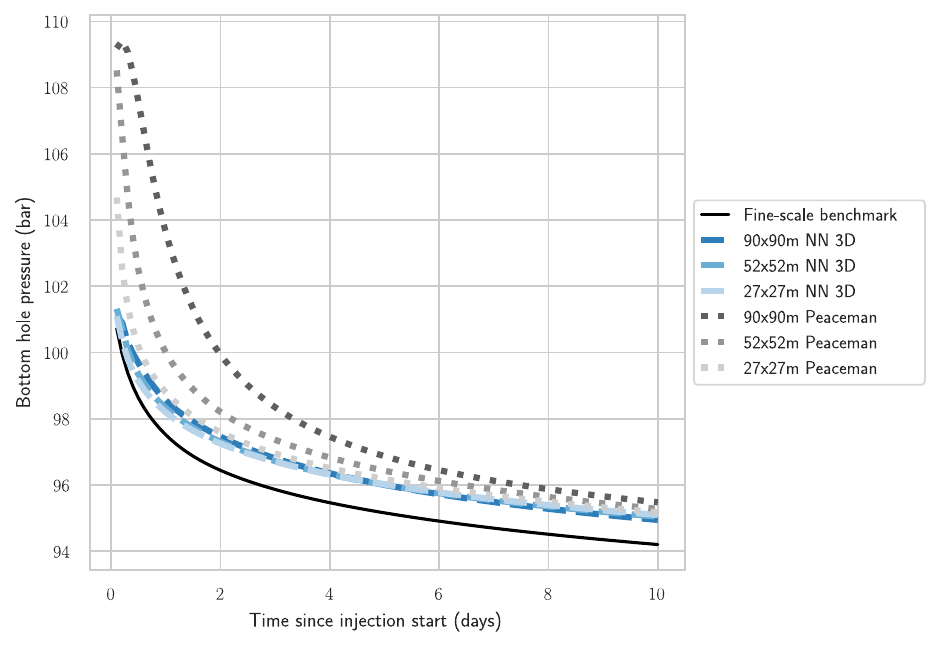}
            \end{minipage}\hfill
            \begin{minipage}[c]{0.3\textwidth}
                \subcaption{Simulation values:
                	\begin{align*}
                		k_{xy,1} &= \qty{8d-12}{\meter\squared} \\
                		k_{xy,2} &= \qty{5d-12}{\meter\squared} \\
                		k_{xy,3} &= \qty{1d-12}{\meter\squared} \\
                		k_{xy,4} &= \qty{8d-13}{\meter\squared} \\
                		k_{xy,5} &= \qty{5d-13}{\meter\squared} \\
                		p_{init} &= \qty{90}{\bar}
                	\end{align*}
                }
                \label{figure:CO2 3D bhp_2}
            \end{minipage}
			\bigskip
            \begin{minipage}[c]{0.67\textwidth}
                \centering
                \includeinkscape[inkscapelatex=false, width=\linewidth]{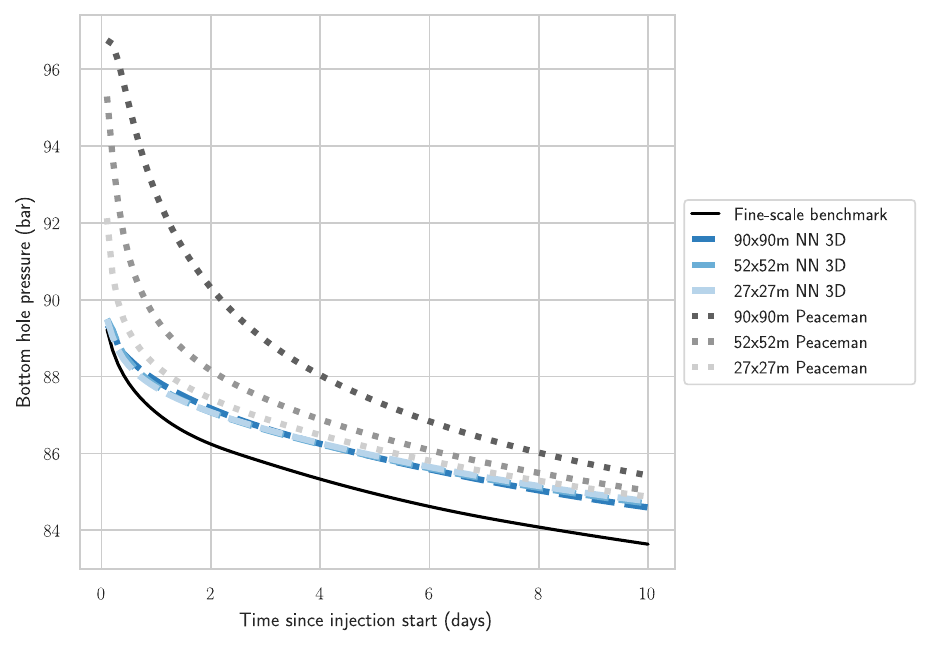}
            \end{minipage}\hfill
            \begin{minipage}[c]{0.3\textwidth}
                \subcaption{Simulation values:
                	\begin{align*}
                		k_{xy,1} &= \qty{9d-12}{\meter\squared} \\
                		k_{xy,2} &= \qty{5d-13}{\meter\squared} \\
                		k_{xy,3} &= \qty{2d-12}{\meter\squared} \\
                		k_{xy,4} &= \qty{5d-12}{\meter\squared} \\
                		k_{xy,5} &= \qty{9d-12}{\meter\squared} \\
                		p_{init} &= \qty{80}{\bar}
                	\end{align*}
                }
                \label{figure:CO2 3D bhp_3}
            \end{minipage}
            \caption[CO\textsubscript{2} 3D \(p_{bhp}\)]{CO\textsubscript{2} 3D example: \(p_{bhp}\) inferred by the ML model and multiphase Peaceman on coarse meshes in comparison with a fine-scale benchmark solution. The ML model is able to infer the fine-scale dynamics much more accurately than the Peaceman model.}
            \label{figure:CO2 3D bhp}
        \end{figure}

\begin{figure}[phbt!]
            \centering
            \begin{minipage}[c]{0.67\textwidth}
                \centering
                \includeinkscape[inkscapelatex=false, width=\linewidth]{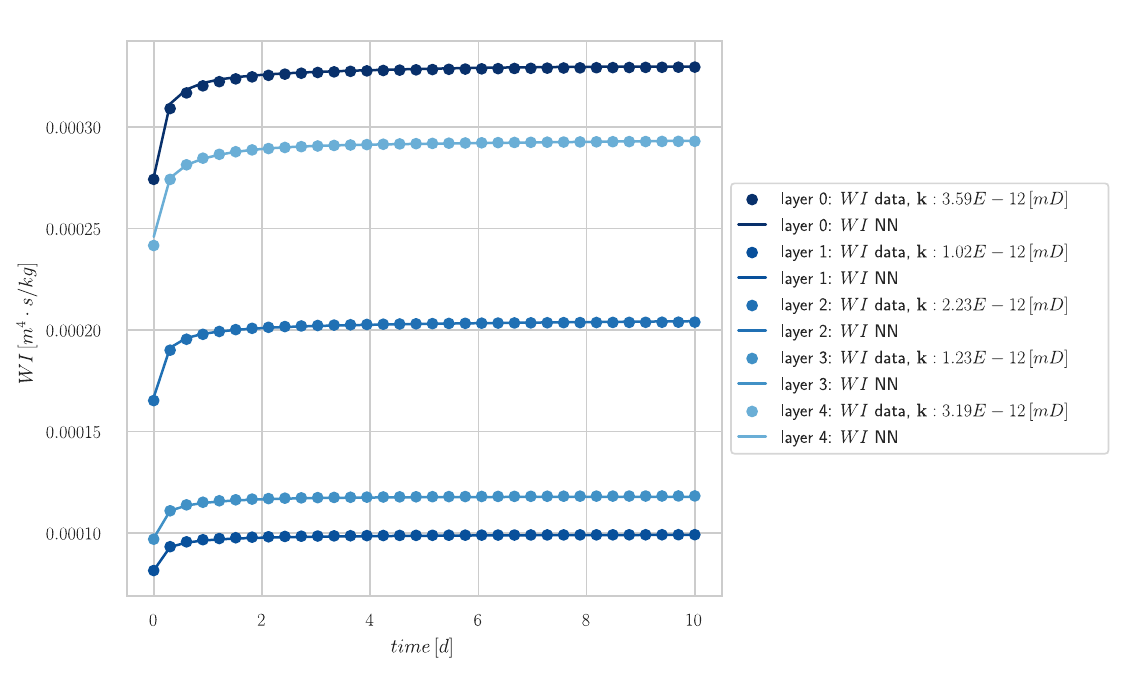}
            \end{minipage}\hfill
            \begin{minipage}[c]{0.3\textwidth}
                \caption{Examplary dependence of \(\widetilde{WI}\) on time for one ensemble member. Both the data points (dots) and the NN prediction (line) are mapped for each of the five layers. It is clearly visible that \(\widetilde{WI}\) increases in the initial days, but remains nearly steady state in the latter two thirds of the injection period.}
                \label{figure:CO2 3D WI vs time}
            \end{minipage}
		\end{figure}

        Input parameters for the ensemble and full-scale simulations are shown in \cref{table:CO2 3D parameters}. We plot the bottom hole pressure over time for three reservoir-scale injection scenarios (\cref{figure:CO2 3D bhp}) with differing geometries and initial pressures, and compare to multiphase Peaceman on different resolutions. On all resolutions, the NN infers the bottom hole pressure accurately and performs significantly better than multiphase Peaceman. The maximum errors are \(\qty{0.9}{\bar},\) \(\qty{1.1}{\bar},\) and \(\qty{1.1}{\bar}\) and the average errors are \(\qty{0.7}{\bar},\) \(\qty{0.9}{\bar},\) and \(\qty{0.9}{\bar}\) for the three simulations, respectively, over the time of injection. In comparison, multiphase Peaceman has maximal errors of \(\qty{23.3}{\bar},\) \(\qty{9.7}{\bar},\) and \(\qty{8.0}{\bar}\) and  average errors of \(\qty{5.4}{\bar},\) \(\qty{2.7}{\bar},\) and \(\qty{3.1}{\bar}\) on the coarsest grid. Only on the finest grid, the results are of comparable accuracy as the ML model. However, even then the Peaceman model infers a too high bottom pressure in the initial days of injection. In the two latter simulation setups (\cref{figure:CO2 3D bhp_2} and \cref{figure:CO2 3D bhp_3}), the error of the ML model increases with time for all grid sizes. It is not clear to the authors what causes this discrepancy, as the model has been trained to model the data accurately (cf. \cref{figure:CO2 3D WI vs time}). Possible causes could be, e.g., that the specific geometries and flow regimes are not as well captured by the training data as the first scenario (\cref{figure:CO2 3D bhp_1}) or that information gets lost in the upscaling process.

        \begin{table}[ht!]
            \centering
            \begin{tabularx}{.9\textwidth}{|>{\raggedright\arraybackslash}X p{1.5cm} >{\raggedright\arraybackslash}X >{\raggedright\arraybackslash}X|}
                \hline
                parameter &  symbol & ensemble values & full simulation values \\
                \hline
                domain size & \(-\)  & \(\qtyproduct{100 x 25}{\meter}\) (radial) &  \(\qtyproduct{1100 x 1100 x 25}{\meter}\) \\
                cell sizes & \(-\)  & log. grid, 50 cells of varying sizes in \(x\)-direction &  \(\qtyproduct{90 x 90 x 5}{\meter}\) (coarse) to \(\qtyproduct{27 x 27 x 5}{\meter}\) (fine) \\
                porosity    & \(\phi\)  & \(\qty{0.35}{}\) & \(\qty{0.35}{}\) \\
                horizontal permeability    & \(k_{xy}\)  & \(\qtyrange[range-units = single, range-phrase =-]{5d-13}{1d-11}{\meter\squared}\) (in 5 layers) & \(\qtyrange[range-units = single, range-phrase =-]{5d-13}{1d-11}{\meter\squared}\) (in 5 layers) \\
                vertical permeability    & \(k_{z}\)  & \(0.5 \cdot k_{xy}\) & \(0.5 \cdot k_{xy}\) \\
                injection rate    & \(Q_t\)  & \(\qty{9342.15}{\tonne\per\day}\) & \(\qty{9342.15}{\tonne\per\day}\) \\
                initial reservoir pressure    & \(p_{init}\)  & \(\qtyrange[range-units = single, range-phrase =-]{50}{120}{\bar}\) & \(\qtyrange[range-units = single, range-phrase =-]{65}{90}{\bar}\) \\
                reservoir temperature & \(T\)  & \(\qty{40}{\degreeCelsius}\) & \(\qty{40}{\degreeCelsius}\) \\
                wellbore radius & \(r_w\)  & \(\qty{0.25}{\meter}\) & \(\qty{0.25}{\meter}\) \\
                \hline
            \end{tabularx}
           	\captionsetup{width=.9\textwidth}
            \caption[CO\textsubscript{2} 3D parameters]{Simulation parameters for the 3D CO\textsubscript{2} example. Note that both the ensemble domain and the full domain have large pore-volume cells at the horizontal boundaries. For the fine-scale benchmark, cell sizes are \(\qtyproduct{4.5 x 4.5 x 1}{\meter}\) in well vicinity and \(\qtyproduct{18 x 18 x 1}{\meter}\) further out.}
            \label{table:CO2 3D parameters}
		\end{table}

\begin{figure}[phbt!]
            \centering
            \includeinkscape[inkscapelatex=false, width=\linewidth]{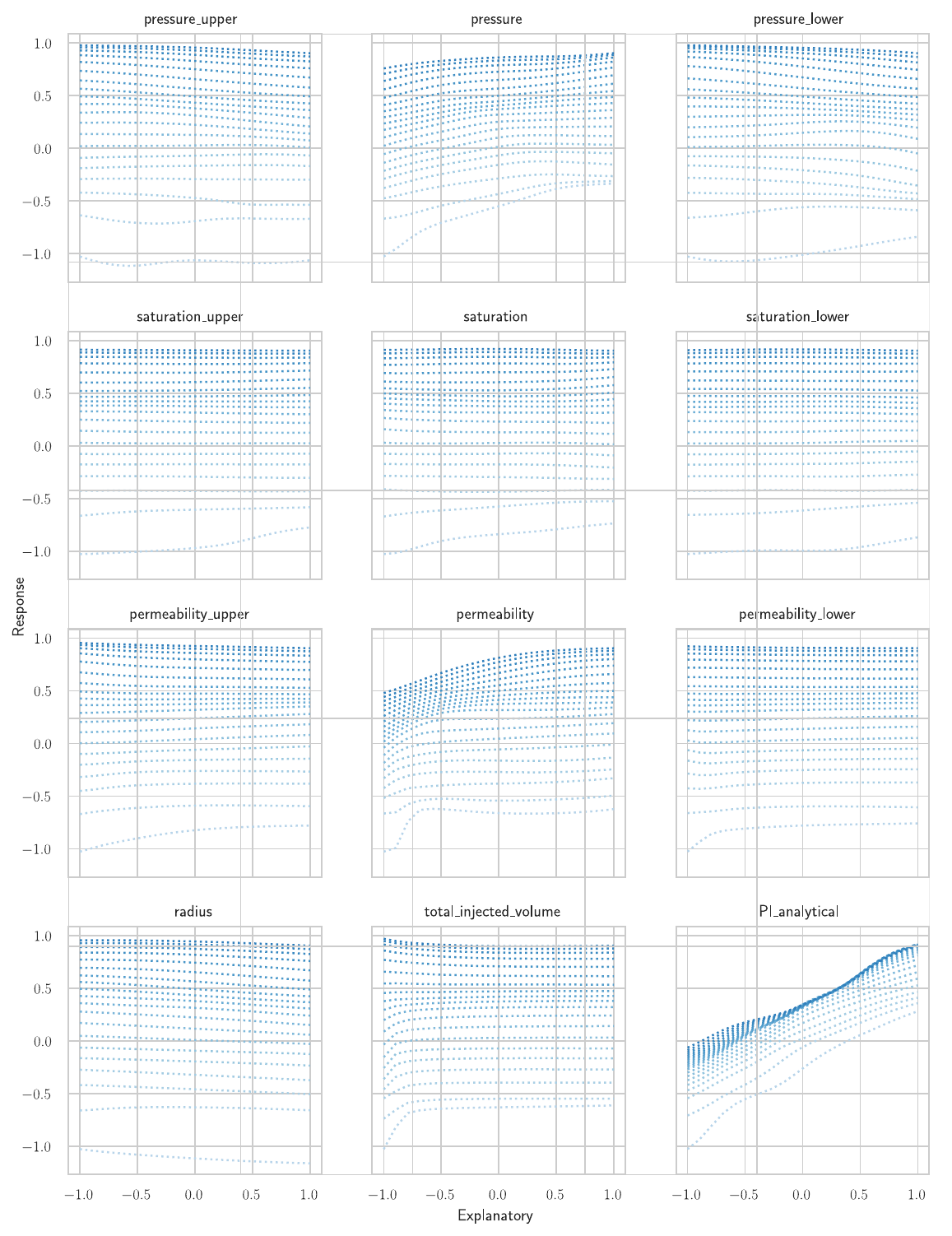}
            \caption{Sensitivity of the model to its inputs. The scaled output of the model is plotted as the denoted scaled input parameter is varied between \(-1\) and \(1.\) All other input variables are fixed to random values. The dotted lines of different shades of blue represent different draws for these random values.}
            \label{figure:CO2 3D sensitivity analysis}
        \end{figure}

        To understand how the model incorporates both the expert knowledge \eqref{equation:expert knowledge} and the correction factor \eqref{equation:correction factor} stemming from phase-mobilities and vertical flow, we study its reliance on the different inputs. \cref{figure:CO2 3D sensitivity analysis} displays the sensitivity of the NN w.r.t. each of its twelve inputs. For all tested combinations of input values, the model relies most on the expert knowledge in terms of the geometrical part of the well index \eqref{equation:expert knowledge}. In comparison, the equivalent well radius, which we know is of importance, has negligible influence. This indicates, that instead of taking the equivalent well radius into account, the model learns the corresponding information from \eqref{equation:expert knowledge}. This displays the advantage of including expert knowledge, as the model does not need to learn the \(\frac{1}{\ln(r_e/r_w)}\) relationship. In addition, the model is sensitive to \(V_{tot},k_{xy},\) and \(p,\) which indicates that they form relevant information the model needs to infer the correction factor \eqref{equation:correction factor}. The total injected volume has an effect on the model's output only for low input values. This is consistent with the fact that the data-driven well index changes fast in the initial three days of injection and remains near constant afterward, as shown in \cref{figure:CO2 3D WI vs time}. This initial period corresponds to the lower third of the input range of \(V_{tot}.\) The model is  slightly sensible to \(p_u\) and \(p_l,\) indicating that it learns the influence of vertical flow on the well injectivity from these two values. The model is nearly insensible to \(k_{xy,u},k_{xy,l},S,S_u,\) and \(S_l,\) indicating that either their influence on the well injectivity is negligible or that the model learns the influence of vertical flow from other inputs. We note that this analysis does not give a comprehensive description of the model's behavior. In particular, the studied combinations of inputs are random and do not necessarily correspond to physical scenarios. Nevertheless, the importance of different inputs in the model's decision-making is indicated.

	\subsection{Discussion of the numerical results}\label{subsection:results discussion}
		The numerical results show that the NN near-well model can learn the Peaceman well model from data in single-phase flow regimes, as well as its ability to infer complex relations in transient two-phase flow regimes. 

		In transient two-phase flow regimes, the model showcases vastly improved accuracy when compared to multiphase Peaceman. It accurately captures fine-scale effects and is able to infer them from coarse-scale simulation values. This leads to increased performance, as the simulation can run on orders of magnitude coarser grids while retaining the same accuracy.

		The model remains flexible, with its exact design being up to the user. Expert knowledge can either be learned from data or explicitly included as a model input. Furthermore, the model can take any subset of simulation values as input, which enables modeling of complex relationships that span multiple coarse grid cells without the need for complex analytic derivation. Studying the final model allows to gain a simplified understanding of the relationships inferred by the model.
		
		Integrated with the OPM Flow-NN framework, the NN near-well model did not lead to notable increase in computation time when compared to multiphase Peaceman.

	\section{Conclusions}\label{section:conclusion}
    	In this work, we show the first integration of NNs into the reservoir simulator OPM Flow. NNs trained in Python within Keras are deployed to OPM Flow by saving the models in an appropriate format and using an automated tool to handle the model interpretation, layer conversion, optimization, and code generation steps. The generated C++ code is then further customized and incorporated into OPM Flow. This provides a seamless transition from the NNs trained within Keras to OPM Flow without requiring significant changes to existing scripts inside OPM Flow. Furthermore, the NNs are fully integrated into OPM's Flow automatic differentiation framework, which allows for effective numerical solution.

		The flexibility of the framework allows for the deployment of NNs in all components of OPM Flow. This facilitates the future development of hybrid models that replace selected parts of a reservoir simulator with ML models to achieve faster simulations and/or higher fidelity to real-world data.

    	The proof-of-concept for a machine-learned near-well model highlights the capabilities of the OPM Flow-NN framework. Trained on a data-driven inflow-performance relation obtained from fine-scale ensemble simulations, the model significantly improves accuracy of the inflow-performance relation in comparison to multiphase Peaceman, while retaining fast inference. The model's flexibility is shown exemplary on two CO\textsubscript{2} injection scenarios that display significant transient effects, which are otherwise only captured by fine-scale resolution of near-well regions. Future extensions of this prototype may include exploring different network architectures such as convolutional neural networks (CNNs) or recurrent neural networks (RNNs) that are better suited for dealing with spatial or temporal relations in complex settings.

	\section*{Acknowledgements}
    	The authors acknowledge funding from the Centre of Sustainable Subsurface Resources (CSSR), grant nr. 331841, supported by the Research Council of Norway, research partners NORCE Norwegian Research Centre and the University of Bergen, and user partners Equinor ASA, Wintershall Dea Norge AS, Sumitomo Corporation, Earth Science Analytics, GCE Ocean Technology, and SLB Scandinavia.

		The first author, Peter von Schultzendorff, acknowledges funding from L. Meltzers Høyskolefond as well as from the Klima og energiomstilling fond of the University of Bergen.

	\section*{Data availability}
		The OPM Flow-NN framework will be published in a future version of OPM Flow. Scripts to reproduce examples for the near-well model are found at the \href{https://github.com/cssr-tools/ML_near_well}{\emph{cssr-tools/ML\_near\_well}} GitHub repository. They require the \href{https://github.com/cssr-tools/pyopmnearwell}{\emph{cssr-tools/pyopmnearwell}} python bindings for near-well modeling.

	\bibliography{references}

\end{document}